\pdfoutput=1
\documentclass[a4paper,12pt]{article} 
\usepackage[top = 2.5cm, bottom = 2.5cm, left = 2.5cm, right = 2.5cm]{geometry}
\usepackage{setspace}
\usepackage{Preamble}

\usepackage{JamesKnot}
\colorlet{basecolor}{red!10!purple!10}
\definecolor{goldaccent}{RGB}{200, 146, 42}

\newcommand{\stdwidth}{1pt}

\newcommand{\setcirclecoords}[1]{
    \coordinate (N) at (0,#1);
    \coordinate (S) at (0,-#1);
    \coordinate (W) at (-#1,0);
    \coordinate (E) at (#1,0);
    \coordinate (NE) at (45:#1);
    \coordinate (NW) at (135:#1);
    \coordinate (SE) at (-45:#1);
    \coordinate (SW) at (-135:#1);
    \coordinate (C) at (0,0);
    \coordinate (CN) at (0,{#1/2});
    \coordinate (CS) at (0,{-#1/2});
    \coordinate (CW) at ({-#1/2},0);
    \coordinate (CE) at ({#1/2},0);
}

\newcommand{\setsquarecoords}[1]{
    \coordinate (N) at (0,#1);
    \coordinate (S) at (0,-#1);
    \coordinate (W) at (-#1,0);
    \coordinate (E) at (#1,0);
    \coordinate (NE) at (#1,#1);
    \coordinate (NW) at (-#1,#1);
    \coordinate (SE) at (#1,-#1);
    \coordinate (SW) at (-#1,-#1);
    \coordinate (C) at (0,0);
    \coordinate (CN) at (0,{#1/2});
    \coordinate (CS) at (0,{-#1/2});
    \coordinate (CW) at ({-#1/2},0);
    \coordinate (CE) at ({#1/2},0);
    \coordinate (CNE) at ({#1/2},{#1/2});
    \coordinate (CNW) at ({-#1/2},{#1/2});
    \coordinate (CSW) at ({-#1/2},{-#1/2});
    \coordinate (CSE) at ({#1/2},{-#1/2});
}

\newcommand{\showgridpoints}[1]{%
  \begin{scope}
    \foreach \x in {-#1,...,#1} {
      \foreach \y in {-#1,...,#1} {
        \pgfmathtruncatemacro{\xeven}{mod(\x,2)==0 ? 1 : 0}
        \pgfmathtruncatemacro{\yeven}{mod(\y,2)==0 ? 1 : 0}
        \ifnum\xeven=1
          \ifnum\yeven=1
            \fill[gray!70] (\x,\y) circle (0.25);
          \fi
        \fi
      }
    }
  \end{scope}
}
\tikzset{
    ->-/.style 2 args={
        postaction={decorate},
        decoration={markings, mark=at position #1 with {\arrow[thick, #2]{>}}}
    },
    ->-/.default={0.5}{}
}

\tikzset{
    oarc/.style={
        basecolor, double=black, double distance=0.5*\stdwidth, line width={2*\stdwidth}
    }
}

\tikzset{
    qarc/.style={
        red!5, double=black, double distance=0.5*\stdwidth, line width={1.5*\stdwidth}
    }
}

\tikzset{
    warc/.style={
        white, double=black, double distance=0.5*\stdwidth, line width={1.5*\stdwidth}
    }
}

\tikzset{
    garc/.style={
        goldaccent!12!white, double=black, double distance=0.5*\stdwidth, line width={1.5*\stdwidth}
    }
}

\tikzset{
    tarc/.style={
        white, double=black, double distance=0.5*\stdwidth, line width={0.75*\stdwidth}
    }
}

\tikzset{
	uarc/.style={
		black, line width=\stdwidth
	}
}

\tikzset{
    uwarc/.style={
        line width=\stdwidth, black, preaction={draw, line width=2*\stdwidth}
    }
}

\tikzset{
    owarc/.style={
        line width=1.75*\stdwidth, white, preaction={preaction={draw, line width=5*\stdwidth, white},draw, line width=2.8pt}
    }
}

\tikzset{
    -<-/.style 2 args={
        postaction={decorate},
        decoration={markings, mark=at position #1 with {\arrow[thick, #2]{<}}} 
    },
    -<-/.default={0.5}{}
}

\tikzset{
    uarcred/.style={
        red, ultra thick
    }
}

\tikzset{
    oarcred/.style={
        basecolor, double=red, double distance=1.6pt, line width=2.4pt
    }
}

\tikzset{
    uwarcred/.style={
        line width=1.6pt, red!15, preaction={draw, line width=2.8pt, red}
    }
}

\tikzset{
    owarcred/.style={
        line width=1.6pt, red!15, preaction={preaction={draw, line width=6.8pt, basecolor},draw, line width=2.8pt, red}
        }
}

\tikzset{
  vJW/.style={
    draw,
    fill=white,
    minimum width=2cm,      
    minimum height=8cm,     
    inner sep=0pt,
    anchor=center
  },
  hJW/.style={
    draw,
    fill=white,
    minimum width=8cm,
    minimum height=2cm,
    inner sep=0pt,
    anchor=center
  }
}

\begin{document}

\title{Twist Formulas for $(m,n)$-Colored Skeins and the Colored Jones Polynomial of Odd 2-Bridge Links}
\author{James Collier\\ \small Yale University \\ \small \texttt{james.collier@aya.yale.edu}}
\date{\today}
\maketitle

\begin{abstract}
 The colored Jones polynomial is a powerful knot invariant which can be calculated diagrammatically using the Kauffman bracket of a knot diagram colored by $n$ (or equivalently, cabled by the Chebyshev polynomial of the second kind $e_n$). Some families of knots, such as 2-bridge links, are created by composing sequences of full and half twists in knot diagrams, and these have formulas for evaluation under the Kauffman bracket. This paper generalizes single-color twist formulas from Yuasa and Masbaum to 2-colored twist formulas for $p$ full twists on 2 strands colored $m$ and $n$. This generalized 2-color twist formula is applied to $(m,n)$-colored odd 2-bridge links $[2a_1, \cdots, 2a_l]$ with $l$ odd, following Yuasa's method, to compute the 2-colored Jones polynomial for this family of links, which reduces to Yuasa's result for single-color 2-bridge links for $m=n$.
\end{abstract}


\tableofcontents
\clearpage
\section{Introduction}
    An important family of knot invariants are the colored Jones polynomials $J_n(K,A)$. These are known to be computable diagrammatically using the Kauffman bracket skein relation. Kauffman brackets are a useful way of turning knot diagrams into polynomials by resolving crossings and other local tangle diagrams.
    
    Diagrammatic skein theory plays an important role in the story of knot invariants. Many of these formulas have an algebraic origin, coming from representations of the quantum group $U_q(\spl_2)$. However, skein theory offers another way to express such formulas for knot diagrams, and such formulas may offer well-behaved asymptotic forms, useful in the case of problems like the volume.
    
    Another invariant appears as the normalized form (by the unknot $O$) of the colored Jones polynomial, and is denoted by $V_{n+1}(K,q)$: $$V_{n+1}(K, q) = \frac{J_{n}(K,q^{1/4})}{J_{n}(O, q^{1/4})}$$
    This normalized invariant is at the heart of the Kashaev-Murakami-Murakami volume conjecture, which states that for a hyperbolic knot $K \subset S^3$ (a knot whose knot complement $S^3 \setminus K$ has a hyperbolic structure), the following is true:
    $$ \underset{n \rightarrow \infty} \lim \frac{2\pi}{n} \log \left| V_n(K, e^{2\pi i/n})\right| = \vol(S^3 \setminus K)$$
    
    In general though, it is not known whether the limit on the left side of the equation exists, and a current topic of research is relating these two invariants for more classes of hyperbolic knots. Currently, the volume conjecture is proven for only  certain families of knots such as hyperbolic twist knots \cite{ChenZhu}, but the conjecture has also been extended to other classes of knots via simplicial volume \cite{MM}, and been solved for torus knots \cite{KT}. 
    
    In general, Habiro showed by studying quantum $\spl_2$ invariants, for $K$ a 0-framed knot and $n \in \zz_{\geq 0}$ that there is a closed form in general for $V_{n+1}(K,q)$ found in \cite[Thm.~3.1]{habiro}.
    $$V_{n+1}(K,q) = \sum_{i =0}^n a_i(K) \prod_{\substack{n+1 - i \ \leq \ j \ \leq \ n+1 +i\\ j \neq n+1}} \{j\}$$

    Note that the term $a_i(K)$ does not depend on the color $n$. In the case of the figure eight knot, it turns out that $a_i(K) = 1$ for all $i \in \zz_{\geq 0}$. There are only a few classes of knots whose explicit $a_i(K)$ has been expressed in this form. Expressing the invariant in this cyclotomic form is desirable, since it lends itself helpful in studying the expression's asymptotics in the volume conjecture (see \cite[Prop.~2.2]{ChenZhu} for the cyclotomic expansion of the colored Jones polynomial of twist knots, which is the starting point for the asymptotic expansion in \cite[Thm.~1.2]{ChenZhu}). When evaluating $q=\exp{(2\pi i/N)}$, we have $\{j\} = 2i \sin(\pi j/N)$, a product of sines which behaves quite nicely. While it is not a focus of this paper to explore this formula, it might be natural to ask what other properties or generalizations of $V_{n+1}(K,q)$ can be explored. A similar notion exists for multi-colored Jones polynomials, and the two color case $J_{m,n}$ indeed has some relevance to the volume conjecture for incomplete hyperbolic metrics.
    
   There are existing formulas in the literature (\cite{yuasa} and \cite{masbaum}) for $p$ full or half twists of $n-$colored Jones-Wenzl idempotents, but not that we had found for two differently colored strands $m,n$. This paper derives a formula for $p$ full twists of $(m,n)$-colored strands using methods found in Yuasa's and Masbaum's papers, and applies it to the family of odd 2-bridge link. This twist formula and its application will be shown at the end of the paper after the introductory build up through basic skein theory and the Jones-Wenzl idempotents.

\section{Basic Knot Theory}
    \textbf{Definition 2.1 $\llbracket$Knots and Diagrams$\rrbracket$:} \ A \textit{knot} is a smooth embedding $K: S^1 \ra S^3$, which we often denote with $K$, referring to the image of the map inside of $S^3$. The corresponding \textit{knot diagram} $D$ is the projection of $K$ to a plane (or another space when specified) with at most double intersections with crossing data. There is a similar notion for \textit{links} which are smooth embeddings $L: \coprod S^1 \ra S^3$.\
    
    \hspace{-0.25in}\textbf{Definition 2.2 $\llbracket$Tangles and Equivalence$\rrbracket$:} \ We often talk about \textit{tangles}, which are local diagrams of knots or links (and may be thought of as smooth embeddings of arcs/circles into a disc with the ends of each arc fixed). Knots, links, and tangles are said to be equivalent if they have an \textit{ambient isotopy} between them, relative to their boundary. Equivalence in this respect is equivalent to having knot/link diagrams being related under a finite sequence of Reidemeister moves I, II, and III, which are shown below as a set of tangle relations:

    $$\textbf{Reidemeister Move I:} \quad \mbox{
        \tikz[baseline=-.6ex, scale=.06]{
            \setcirclecoords{12}
            \fill[basecolor] (C) circle [radius=12cm];
            \begin{scope}
                \clip (C) circle [radius=12cm];
                \draw[oarc] (N) to[out=-90, in=west] (3,-2);
                \draw[oarc] (3,-2) to[out=east, in=east, looseness=1.5] (3, 2);
                \draw[oarc] (3,2) to[out=west, in=north] (S);
            \end{scope}
            \draw[dashed] (C) circle [radius=12cm];
        }
    } \leftrightarrow
    \mbox{
        \tikz[baseline=-.6ex, scale=.06]{
            \setcirclecoords{12}
            \fill[basecolor] (C) circle [radius=12cm];
            \begin{scope}
                \clip (C) circle [radius=12cm];
                \draw[oarc] (S) -- (N);
            \end{scope}
            \draw[dashed] (C) circle [radius=12cm];
        }
    } \leftrightarrow
    \mbox{
        \tikz[baseline=-.6ex, scale=.06]{
            \setcirclecoords{12}
            \fill[basecolor] (C) circle [radius=12cm];
            \begin{scope}
                \clip (C) circle [radius=12cm];
                \draw[oarc] (S) to[out=north, in=west] (3,2);
                \draw[oarc] (3,2) to[out=east, in=east, looseness=1.5] (3, -2);
                \draw[oarc] (3,-2) to[out=west, in=south] (N);
            \end{scope}
            \draw[dashed] (C) circle [radius=12cm];
        }
    } 
    $$
    $$\textbf{Reidemeister Move II:} \quad \mbox{
        \tikz[baseline=-.6ex, scale=.06]{
            \setcirclecoords{12}
            \fill[basecolor] (C) circle [radius=12cm];
            \begin{scope}
                \clip (C) circle [radius=12cm];
                \draw[oarc] (NW) to[out=-30, in=90] (3,0) to[out=-90, in=30] (SW);
                \draw[oarc] (NE) to[out=-150, in=90] (-3,0) to[out=-90, in=150] (SE);
            \end{scope}
            \draw[dashed] (C) circle [radius=12cm];
        }
    } \leftrightarrow
    \mbox{
        \tikz[baseline=-.6ex, scale=.06]{
            \setcirclecoords{12}
            \fill[basecolor] (C) circle [radius=12cm];
            \begin{scope}
                \clip (C) circle [radius=12cm];
                \draw[oarc] (NW) to[out=-45, in=45] (SW);
                \draw[oarc] (NE) to[out=-135, in=135] (SE);
            \end{scope}
            \draw[dashed] (C) circle [radius=12cm];
        }
    } \leftrightarrow
    \mbox{
        \tikz[baseline=-.6ex, scale=.06]{
            \setcirclecoords{12}
            \fill[basecolor] (C) circle [radius=12cm];
            \begin{scope}
                \clip (C) circle [radius=12cm];
                \draw[oarc] (NE) to[out=-150, in=90] (-3,0) to[out=-90, in=150] (SE);
                \draw[oarc] (NW) to[out=-30, in=90] (3,0) to[out=-90, in=30] (SW);
            \end{scope}
            \draw[dashed] (C) circle [radius=12cm];
        }
    } 
    $$
    $$\textbf{Reidemeister Move III:} \quad \mbox{
        \tikz[baseline=-.6ex, scale=.06]{
            \setcirclecoords{12}
            \fill[basecolor] (C) circle [radius=12cm];
            \begin{scope}
                \clip (C) circle [radius=12cm];
                \draw[oarc] (W) to[out=45, in=135] (E);
                \draw[oarc] (SW) -- (NE);
                \draw[oarc] (SE) -- (NW);
            \end{scope}
            \draw[dashed] (C) circle [radius=12cm];
        }
    } \leftrightarrow
    \mbox{
        \tikz[baseline=-.6ex, scale=.06]{
            \setcirclecoords{12}
            \fill[basecolor] (C) circle [radius=12cm];
            \begin{scope}
                \clip (C) circle [radius=12cm];
                \draw[oarc] (W) to[out=-45, in=-135] (E);
                \draw[oarc] (SW) -- (NE);
                \draw[oarc] (SE) -- (NW);
            \end{scope}
            \draw[dashed] (C) circle [radius=12cm];
        }
    }$$

    Links and knots which are only equivalent under Reidemeister move II and III are said to be \textit{regularly isotopic}. Knot invariants are functions from the space of knots to another space (such as to Laurent polynomials e.g. $\zz[A^{\pm1}]$) which preserve equivalence under ambient isotopy. Invariants often allow for classification and the distinguishing of complicated knots or knots with large crossing number. A function on knot diagrams is invariant when it remains unchanged under R I-III moves on the diagram.\

    \hspace{-0.25in}\textbf{Definition 2.3 $\llbracket$$q$-Calculus Notation$\rrbracket$:} \ 
    Here we will describe some of the commonly used notations for quantum integers using $k,n \in \zz$ and $q$ an indeterminate. However, it should be noted that some authors prefer different variables. In this paper we will use $q$ with $q = a^2 = A^4$.
    $$\text{With }\{k\} = q^{k/2} - q^{-k/2} \text{ we define the q-bracket: }[k] = \frac{\{k\}}{\{1\}}$$
    Frequently, we will find the need to use an analogue of factorials for the $q$ bracket.
    $$\{k\}! = \{k\}\{k-1\} \cdots \{1\} \mathand [k]! = [k][k-1]\cdots [1]$$
    The $q$-bracket has its associated $q$-binomial coefficient:
    $$\bracketbinom{n}{k}_q = \frac{[n]!}{[k]![n-k]!}$$
    We also define the q-Pochhammer symbol and its q-binomial coefficients, which are said to be q-analogues of their original versions:
    $$(q; q)_k = (q)_k \coloneq \prod_{i = 1}^k(1-q^i) \mathand \binom{n}{k}_q = \frac{(q)_n}{(q)_k (q)_{n-k}}$$
    There are multinomial forms as well:
    $$\binom{n}{k_1, k_2, \cdots, k_m}_q = \frac{(q)_n}{(q)_{k_1}(q)_{k_2}\cdots (q)_{k_m}} \qquad \bracketbinom{n}{k_1, k_2, \cdots, k_m}_q = \frac{[n]!}{[k_1]![k_2]!\cdots[k_m]!}$$
    A defining motivation behind q-calculus and its functions are that when the indeterminate $q \ra 1$, the functions reduce to their classical analogues. These numbers arise in the representations of $\spl_2$ with its quantum group $U_q \spl_2$, but it is beyond the scope of this paper, which considers the graphical aspects associated with skein relations and knot invariants.

\section{The Kauffman Bracket Skein Module}
    \textbf{Definition 3.1 $\llbracket$Web Spaces$\rrbracket$:} \ 
    Denote $T_{2n}$ the space of tangle diagrams in $D^2$, the unit disk with $n$ fixed arcs (having $2n$ distinct endpoints lying on the boundary $\partial D^2 = S^1$), and $W_{2n}$ the free $\zz[q^{\pm 1/4}]$-module with basis consisting of non-intersecting tangle diagrams in $T_{2n}$ up to ambient isotopy. Here $\zz[q^{\pm 1/4}] = \zz[q^{1/4}, q^{-1/4}]$ is the ring of Laurent polynomials in the indeterminate $q^{1/4}$ with integer coefficients. $W_{2n}$ is called a \textit{web space}.\

    \hspace{-0.25in}\textbf{Definition 3.2 $\llbracket$Kauffman Bracket$\rrbracket$}
    Then the \textit{Kauffman Bracket} is defined as a $\zz[q^{\pm 1/4}]$-linear map $\langle \, \cdot \,  \rangle:   \zz[q^{\pm 1/4}] T_{2n} \ra  W_{2n}$ with the following properties for local tangle diagrams:
    $$(1) \qquad \kbrack{\mbox{
        \tikz[baseline=-.6ex, scale=.06]{
            \setcirclecoords{12}
            \fill[basecolor] (C) circle [radius=12cm];
            \begin{scope}
                \clip (C) circle [radius=12cm];
                \draw[oarc] (NW) -- (SE);
                \draw[oarc] (SW) -- (NE);
            \end{scope}
            \draw[dashed] (C) circle [radius=12cm];
        }
    }} = q^{1/4} \kbrack{\mbox{
        \tikz[baseline=-.6ex, scale=.06]{
            \setcirclecoords{12}
            \fill[basecolor] (C) circle [radius=12cm];
            \begin{scope}
                \clip (C) circle [radius=12cm];
                \draw[oarc] (NW) to[out=-45, in=45] (SW);
                \draw[oarc] (NE) to[out=-135, in=135] (SE);
            \end{scope}
            \draw[dashed] (C) circle [radius=12cm];
        }
    }} + q^{-1/4} \kbrack{\mbox{
        \tikz[baseline=-.6ex, scale=.06]{
            \setcirclecoords{12}
            \fill[basecolor] (C) circle [radius=12cm];
            \begin{scope}
                \clip (C) circle [radius=12cm];
                \draw[oarc] (NW) to[out=-45, in=-135] (NE);
                \draw[oarc] (SW) to[out=45, in=135] (SE);
            \end{scope}
            \draw[dashed] (C) circle [radius=12cm];
        }
    }}$$
    $$(2) \qquad \kbrack{D \sqcup \mbox{
        \tikz[baseline=-.6ex, scale=.06]{
            \setcirclecoords{12}
            \fill[basecolor] (C) circle [radius=12cm];
            \begin{scope}
                \clip (C) circle [radius=12cm];
                \draw[oarc] (C) circle [radius=8cm];
            \end{scope}
            \draw[dashed] (C) circle [radius=12cm];
        }
    }} = -[2] \kbrack{D} = (-q^{1/2} - q^{-1/2}) \kbrack{D}$$\
    
    \hspace{-0.25in}\textbf{Definition 3.3: $\llbracket$Kauffman Bracket Skein Module on the Annulus$\rrbracket$}: \ It is also useful to imagine knots and tangles embedded in different surfaces. One particularly useful model is on the annulus $S^1 \times I$. The \textit{Kauffman Bracket Skein Module} (KBSM) of the annulus is denoted $K(S^1 \times I)$, and is the $\zz[q^{\pm1/4}]$-module generated by link diagrams in the annulus, modulo the Kauffman bracket skein relations of Definition 3.2. It is isomorphic to $\zz[q^{\pm1/4}][z]$ where $z$ denotes the embedding of the circle around the hole of the annulus, and $z^n$ is the $n$-cabling.
    $$1 = \varnothing = \mbox{
        \tikz[baseline=-.6ex, scale=.06]{
            \setcirclecoords{12}
            \fill[basecolor] (C) circle [radius=12cm];
            \begin{scope}
                \clip (C) circle [radius=12cm];
                \filldraw[color = red!5, draw=black, dashed] (C) circle [radius=3cm];
            \end{scope}
            \draw[dashed] (C) circle [radius=12cm];
        }
    } \quad z = \mbox{
        \tikz[baseline=-.6ex, scale=.06]{
            \setcirclecoords{12}
            \fill[basecolor] (C) circle [radius=12cm];
            \begin{scope}
                \clip (C) circle [radius=12cm];
                \filldraw[color = red!5, draw=black, dashed] (C) circle [radius=3cm];
                \draw[oarc] (C) circle [radius=7.5cm];
            \end{scope}
            \draw[dashed] (C) circle [radius=12cm];
        }
    } \quad z^3 = \mbox{
        \tikz[baseline=-.6ex, scale=.06]{
            \setcirclecoords{12}
            \fill[basecolor] (C) circle [radius=12cm];
            \begin{scope}
                \clip (C) circle [radius=12cm];
                \filldraw[color = red!5, draw=black, dashed] (C) circle [radius=3cm];
                \draw (C) circle [radius=6.5cm];
                \draw (C) circle [radius=8cm];
                \draw (C) circle [radius=9.5cm];
            \end{scope}
            \draw[dashed] (C) circle [radius=12cm];
        }
    }$$

    An important set of elements arise in $K(S^1 \times I)$ as the eigenvectors of the twist operator $t: K(S^1 \times I) \ra K(S^1 \times I)$ which induces a positive "right-handed" twist on a diagram in $K(S^1 \times I)$, then resolves the induced crossings with the Kauffman bracket. 
    The eigenvectors of the twist operator are actually a family of polynomials defined recursively:
    $$e_0 = 1 \quad e_1 = z \quad e_n = ze_{n-1} - e_{n-2}$$
    These have eigenvalues under the twist operator:
    $$ t (e_n) = \mu_n e_n \quad \text{where } \quad \mu_n = (-1)^n A^{n^2 + 2n} = (-1)^n q^{(n^2 + 2n)/4}$$
    
    These are called the \textit{Chebyshev polynomials of the second kind} (whereas the first kind starts with $e_0 = 2$), and they play an important role in skein theory. A strand is said to be $n-$colored when it is cabled by $e_n$ in $K(S^1 \times I)$. These polynomials are also used to define the colored Jones polynomial, which is an invariant for any $n \in \zz_{\geq 0}$ and knot $K$ with:
    $$J_n(K,A) = \left((-1)^n A^{n^2 + 2n}\right)^{-w(D)}\langle e_n\rangle_D$$
    Where $w(D)$ is the \textit{writhe} of a knot diagram, a weighted sum of positively and negatively oriented crossings, and $\langle e_n\rangle_D$ denotes the $e_n-$cabling of the knot diagram $D$. Both the writhe and the Kauffman bracket are invariant under regular isotopy (Reidemeister moves II and III). The colored Jones polynomials are constructed to be additionally invariant under ambient isotopy by introducing the twist eigenvalue to counterbalance the coefficients induced by Reidemeister move I.

\section{Temperley-Lieb Algebra \& JW Idempotents}
    \textbf{Definition 4.1 $\llbracket$Temperley-Lieb algebra$\rrbracket$:}
    An analogous way to think about colored strands is through the lens of the \textit{Temperley-Lieb algebra} $TL_{n}$. Consider diagrams in a square with $n$ marked points on the left edge (inputs) and $n$ marked points on the right edge (outputs), consisting of $n$ pairwise disjoint arcs joining these $2n$ points in pairs, with no crossings and no closed loops. An arc may join an input to an output, two inputs, or two outputs. $TL_n$ is the free $\zz[q^{\pm 1/4}]$-module with basis the isotopy classes (rel boundary) of such diagrams; there are $\frac{1}{n+1}\binom{2n}{n}$ of them (the $n$th Catalan number). Here is a diagram showing examples of elements in $TL_3$.

    $$I: \mbox{
        \tikz[baseline=-.6ex, scale=.075]{
            \setsquarecoords{12}
             
            \begin{scope}
                 
                \draw[warc] (-12, 8) -- (12,8);
                \draw[warc] (-12, 0) -- (12,0);
                \draw[warc] (-12, -8) -- (12, -8);
                \draw[densely dotted] (-12,-12) -- (-12, 12);
                \draw[densely dotted] (12,-12) -- (12, 12);

            \end{scope}
        }
    }\quad U_1: \mbox{
        \tikz[baseline=-.6ex, scale=.075]{
            \setsquarecoords{12}
             
            \begin{scope}

                \draw[warc] (-12, 8) to[out=east, in= east, looseness=2] (-12,0);
                \draw[warc] (12, 8) to[out=west, in= west, looseness=2] (12,0);
                \draw[warc] (-12, -8) -- (12, -8);
                \draw[densely dotted] (-12,-12) -- (-12, 12);
                \draw[densely dotted] (12,-12) -- (12, 12);

            \end{scope}
        }
    }\quad U_2: \mbox{
        \tikz[baseline=-.6ex, scale=.075]{
            \setsquarecoords{12}
             
            \begin{scope}
                 
                \draw[warc] (-12, 0) to[out=east, in= east, looseness=2] (-12,-8);
                \draw[warc] (12, 0) to[out=west, in= west, looseness=2] (12,-8);
                \draw[warc] (-12, 8) -- (12, 8);
                \draw[densely dotted] (-12,-12) -- (-12, 12);
                \draw[densely dotted] (12,-12) -- (12, 12);

            \end{scope}
        }
    }$$
    $$U_1 \cdot U_2: \mbox{
        \tikz[baseline=-.6ex, scale=.075]{
            \setsquarecoords{12}
             
            \begin{scope}
                 
                \draw[warc] (-12, 8) to[out=east, in= east, looseness=2] (-12,0);
                \draw[warc] (12, 0) to[out=west, in= west, looseness=2] (12,-8);
                \draw[warc] (-12, -8) to[out=east, in=west] (12, 8);
                \draw[densely dotted] (-12,-12) -- (-12, 12);
                \draw[densely dotted] (12,-12) -- (12, 12);

            \end{scope}
        }
    } = \mbox{
        \tikz[baseline=-.6ex, scale=.075]{
            \setsquarecoords{12}
             
            \begin{scope}

                \draw[warc] (-12, 8) to[out=east, in= east, looseness=2] (-12,0);
                \draw[warc] (12, 8) to[out=west, in= west, looseness=2] (12,0);
                \draw[warc] (-12, -8) -- (12, -8);
                \draw[densely dotted] (-12,-12) -- (-12, 12);
                \draw[densely dotted] (12,-12) -- (12, 12);

            \end{scope}
        }
    } \mkern-16.5mu
    \mbox{
        \tikz[baseline=-.6ex, scale=.075]{
            \setsquarecoords{12}
             
            \begin{scope}
                 
                \draw[warc] (-12, 0) to[out=east, in= east, looseness=2] (-12,-8);
                \draw[warc] (12, 0) to[out=west, in= west, looseness=2] (12,-8);
                \draw[warc] (-12, 8) -- (12, 8);
                \draw[densely dotted] (-12,-12) -- (-12, 12);
                \draw[densely dotted] (12,-12) -- (12, 12);

            \end{scope}
        }
    }$$
    Multiplication in $TL_n$ is given on basis diagrams by stacking: the outputs of the first diagram are glued to the inputs of the second, as can be seen above, and each closed loop created is removed and replaced by the factor $-[2]$. This is extended bilinearly. The identity $I = 1_n$ is the diagram of $n$ vertical arcs, and $U_i$ (for $1 \leq i \leq n-1$) is the diagram joining the $i$th and $(i+1)$th inputs to each other, the $i$th and $(i+1)$th outputs to each other, and every other input vertically to its output. As an algebra, $TL_n$ is generated by $1_n, U_1, \cdots, U_{n-1}$, which satisfy
    $$U_i^2 = -[2]\,U_i, \qquad U_iU_{i\pm 1}U_i = U_i, \qquad U_iU_j = U_jU_i \ \text{ for } |i-j| \geq 2,$$
    so every element of $TL_n$ is a linear combination of products of the $U_i$'s and $1_n$. We will switch between the notion of Temperley-Lieb algebra and the web space $W_{2n}$ quite often, as they are similar notions, but $TL_n$ offers an introduction into Jones-Wenzl idempotents.
    
    As with $W_{2n}$, any tangle diagram in the square with $n$ inputs and $n$ outputs, possibly with crossings, determines an element of $TL_n$ by applying the Kauffman bracket relations to resolve every crossing and remove every closed loop.

    \hspace{-0.25in}\textbf{Definition 4.2 $\llbracket$Jones-Wenzl Idempotents$\rrbracket$:} \ We define the \textit{Jones-Wenzl idempotents} $f_k \in TL_k$ recursively by
    $$f_1 = 1_1 \qquad f_{k+1} = f_{k} + \frac{[k]}{[k+1]}f_kU_{k}f_k,$$
    where on the right $f_k$ is regarded as an element of $TL_{k+1}$ by adding one vertical strand on the bottom. Because the recursion divides by quantum integers, the coefficients of $f_k$ are in general rational functions rather than Laurent polynomials (already $f_2 = 1_2 + \frac{1}{[2]}U_1$). So $f_k$ lives in $TL_k \otimes \Q(q^{1/4})$, where $\Q(q^{1/4})$ denotes the field of rational functions in $q^{1/4}$. From here on, whenever Jones-Wenzl idempotents appear, we work over $\Q(q^{1/4})$, that is, in $TL_n \otimes \Q(q^{1/4})$ and $W_{2n} \otimes \Q(q^{1/4})$. The $f_k$ satisfy the following properties \cite{kauffman}:
    
    \begin{enumerate}
        \item \textbf{(Idempotency)} $$f_k^2 = f_k$$
        \item \textbf{(Annihilation of Turnbacks)} $$f_kU_j = U_jf_k = 0$$
        \item \textbf{(Closure Property)} Define $\tr_q:TL_n \otimes \Q(q^{1/4}) \ra \Q(q^{1/4})$ which takes an element of $TL_n$, closes each strand by connecting the $i$'th input and output, then taking the Kauffman bracket of the entire diagram. Then:
        $$\tr_q(f_k) = (-1)^k[k+1],$$
        which is again a Laurent polynomial.
    \end{enumerate}
    
    The Jones-Wenzl idempotent $f_n$ can also be drawn as a white box labeled $n$ on $n$ parallel strands. The recursion above then takes the following graphical form:
$$\kbrack{\mkern-8mu \mbox{
    \tikz[baseline=-.6ex, scale=.075]{
        \setsquarecoords{12}
         
        \begin{scope}
            \draw[warc] (-10, 0) -- (10, 0) node[midway, above, black, yshift=6pt, xshift=0pt]{\scriptsize $n+1$};
            \node[vJW, scale=0.08] at (0,0) {};
        \end{scope}
    }
}\mkern-8mu} = \kbrack{ \mkern-8mu \mbox{
    \tikz[baseline=-.6ex, scale=.075]{
        \setsquarecoords{12}
         
        \begin{scope}
            \draw[warc] (-12, 6) -- (12, 6) node[midway, above, black, yshift=0pt, xshift=7pt]{\scriptsize $n$};
            \node[vJW, scale=0.1] at (0,6) {};
            \draw[warc] (-12, -6) -- (12, -6) node[midway, below, black, yshift=0pt, xshift=0pt]{\scriptsize $1$};
        \end{scope}
    }
}\mkern-8mu} + \frac{[n]}{[n+1]}\kbrack{ \mkern-8mu \mbox{
    \tikz[baseline=-.6ex, scale=.075]{
        \setsquarecoords{12}
         
        \begin{scope}
            \draw[warc] (-15, 4) -- (-8, 4) node[midway, above, black, yshift=0pt, xshift=0pt]{\scriptsize $n$};
            \draw[warc] (8, 4) -- (15, 4) node[midway, above, black, yshift=0pt, xshift=0pt]{\scriptsize $n$};

            \draw[warc] (-15, -6) -- (-8, -6) node[midway, above, black, yshift=-2pt, xshift=0pt]{\scriptsize $1$};
            \draw[warc] (8, -6) -- (15, -6) node[midway, above, black, yshift=-2pt, xshift=0pt]{\scriptsize $1$};
            
            \draw[warc] (-8,2) to[out=east, in=east, looseness=1.5] (-8, -6);
            \draw[warc] (8,2) to[out=west, in=west, looseness=1.5] (8, -6);

            \draw[warc] (-8, 6) -- (8, 6) node[midway, above, black, yshift=-2pt, xshift=0pt]{\scriptsize $n-1$};
            \node[vJW, scale=0.1] at (-8,4) {};
            \node[vJW, scale=0.1] at (8,4) {};
        \end{scope}
    }
}\mkern-8mu}$$
The coloring $1$ is simply a single strand. Each box labeled with $n$ can be seen as a 'black box' with $n$ inputs and $n$ outputs. These boxes are also called \textit{clasps}. Each strand labeled with $n$ can be imagined as containing $n$ separate strands.\

\hspace{-0.25in}\textbf{Proposition 4.3:} \cite{kauffman} The Jones-Wenzl idempotents also have the graphical properties: 
    \begin{enumerate}
        \item $\kbrack{\mkern-8mu \mbox{
    \tikz[baseline=-.6ex, scale=.075]{
        \setsquarecoords{12}
         
        \begin{scope}
            \draw[warc] (-10, 0) -- (10, 0);
            \node[vJW, scale=0.1, label={[xshift=0pt, yshift=-1pt]above:\scriptsize $n$}] at (-4,0) {};
            \node[vJW, scale=0.1, label={[xshift=0pt, yshift=-1pt]above:\scriptsize $n$}] at (4,0) {};
        \end{scope}
    }
}\mkern-8mu} = \mbox{
    \tikz[baseline=-.6ex, scale=.075]{
        \setsquarecoords{12}
         
        \begin{scope}
            \draw[warc] (-10, 0) -- (10, 0);
            \node[vJW, scale=0.1, label={[xshift=0pt, yshift=-1pt]above:\scriptsize $n$}] at (0,0) {};
        \end{scope}
    }
}$
    \item $\kbrack{\mkern-12mu \mbox{
    \tikz[baseline=-.6ex, scale=.075]{
        \setsquarecoords{12}
         
        \begin{scope}
            \draw[warc] (-12, 0) -- (-6, 0);
            \draw[warc] (-6, 4) -- (12, 4);
            \draw[warc] (-5, 2) -- (-3, 2);
            \draw[warc] (-5, -2) -- (-3, -2);
            \draw[warc] (-3, 2) to[out=east, in=east, looseness=2] (-3, -2) ;
            \node at (12,4) [right] {\scriptsize $k$};
             \node at (-2,0) [right] {\scriptsize $1$};
            \node at (12,-4) [right] {\scriptsize $n-k-2$};
            \draw[warc] (-6, -4) -- (12, -4);
            \node[vJW, scale=0.1, label={[xshift=0pt, yshift=-1pt]above:\scriptsize $n$}] at (-6,0) {};
        \end{scope}
    }
}\mkern-12mu} = 0, \quad k \in \{0, \cdots, n-2\}$
    \item $\kbrack{\mkern-12mu \mbox{
    \tikz[baseline=-.6ex, scale=.075]{
        \setsquarecoords{12}
         
        \begin{scope}
            \draw[warc] circle [radius=6cm];
            \node[hJW, scale=0.08, label={[xshift=5pt, yshift=-1pt]above:\scriptsize $n$}] at (6,0) {};
        \end{scope}
    }
}\mkern-12mu} = (-1)^n [n+1]$
    \end{enumerate}

    \hspace{-0.25in}\textbf{Remark 4.4:}
    There is some conflation between cabling with $e_n$ and coloring a strand with the Jones-Wenzl idempotent with the label $n$. The convention that we will take is that for a knot K colored with the Jones-Wenzl idempotent $n$, then:
    $$\kbrack{f_n} = \langle e_{n}\rangle $$
    That is, coloring a strand with the label $n$ is the same as cabling a diagram with $e_{n}$.\

    \hspace{-0.25in}\textbf{Proposition 4.5:} \cite{kauffman} We have the following properties which help us understand how to manipulate tangle diagrams containing JW idempotents. Namely, we are able to 'split' JW clasps into smaller clasps and combine them in the same way using (1), resolve self-crossings of strands using (2), resolve self-twists using (3), and resolve partial closures using (4). For $k = 0, 1, \cdots, n$, we have the following:
    \begin{enumerate}
        \item $\kbrack{\mkern-12mu\mbox{
    \tikz[baseline=-.6ex, scale=.075]{
        \setsquarecoords{12}
         
        \begin{scope}
            \draw[warc] (-12, 0) -- (-4, 0);
            \draw[warc] (-4, 4) -- (10, 4);
            \draw[warc] (-4, -4) -- (10, -4);
            \node[vJW, scale=0.1, label={[xshift=0pt, yshift=-1pt]above:\scriptsize $n$}] at (-4,0) {};
            \node at (9,4) [right] {\scriptsize $k$};
            \node at (9,-4) [right] {\scriptsize $n-k$};
            
        \end{scope}
    }
}\mkern-12mu} = \kbrack{\mkern-12mu\mbox{
    \tikz[baseline=-.6ex, scale=.075]{
        \setsquarecoords{12}
         
        \begin{scope}
            \draw[warc] (-12, 0) -- (-4, 0);
            \draw[warc] (-4, 4) -- (10, 4);
            \draw[warc] (-4, -4) -- (10, -4);
            \node[vJW, scale=0.1, label={[xshift=0pt, yshift=-1pt]above:\scriptsize $n$}] at (-4,0) {};
            \node[vJW, scale=0.06, label={[xshift=0pt, yshift=-1pt]above:\scriptsize $k$}] at (4,4) {};
            \node at (9,4) [right] {\scriptsize $k$};
            \node at (9,-4) [right] {\scriptsize $n-k$};
        \end{scope}
    }
}\mkern-12mu}$
    \item $\kbrack{\mkern-12mu\mbox{
    \tikz[baseline=-.6ex, scale=.075]{
        \setsquarecoords{12}
         
        \begin{scope}
            \draw[warc] (-12, 0) -- (-4, 0);
            \draw[warc] (-4, -4) to[out=east, in=west] (10, 4);
            \draw[warc] (-4, 4) to[out=east, in=west] (10, -4);
            \node[vJW, scale=0.1, label={[xshift=0pt, yshift=-1pt]above:\scriptsize $n$}] at (-4,0) {};
            \node at (9,4) [right] {\scriptsize $k$};
            \node at (9,-4) [right] {\scriptsize $n-k$};
            
        \end{scope}
    }
}\mkern-12mu} = q^{k(n-k)/4}\kbrack{\mkern-12mu\mbox{
    \tikz[baseline=-.6ex, scale=.075]{
        \setsquarecoords{12}
         
        \begin{scope}
            \draw[warc] (-12, 0) -- (-4, 0);
            \draw[warc] (-4, 0) -- (4, 0);
            \node[vJW, scale=0.1, label={[xshift=0pt, yshift=-1pt]above:\scriptsize $n$}] at (-4,0) {};
        \end{scope}
    }
}\mkern-12mu}$
    \item $\kbrack{\mkern-12mu\mbox{
    \tikz[baseline=-.6ex, scale=.075]{
        \setsquarecoords{12}
         
        \begin{scope}
            \draw[warc] (-12, 0) -- (-4, 0);
            \draw[warc] (-4, 0) -- (2, 0); 
            \draw[warc] (2, 0) to[out=east, in=east, looseness=1.5] (4, -4);
            \draw[warc] (4, -4) to[out=west, in=west, looseness=1.5] (6, 0);
            \draw[warc] (6, 0) -- (10, 0);
            \node[vJW, scale=0.1, label={[xshift=0pt, yshift=-1pt]above:\scriptsize $n$}] at (-4,0) {};
            \node at (9,0) [right] {\scriptsize $n$};
        \end{scope}
    }
}\mkern-12mu} = \mu_n \kbrack{\mkern-12mu\mbox{
    \tikz[baseline=-.6ex, scale=.075]{
        \setsquarecoords{12}
         
        \begin{scope}
            \draw[warc] (-12, 0) -- (-4, 0);
            \draw[warc] (-4, 0) -- (4, 0);
            \node[vJW, scale=0.1, label={[xshift=0pt, yshift=-1pt]above:\scriptsize $n$}] at (-4,0) {};
        \end{scope}
    }
}\mkern-12mu} \qquad \text{where $\mu_n = (-1)^n q^{(n^2 + 2n)/4}$}$
    \item $\kbrack{\mkern-12mu\mbox{
    \tikz[baseline=-.6ex, scale=.075]{
        \setsquarecoords{12}
         
        \begin{scope}
            \draw[warc] (-12, 2) -- (12, 2); 
            \draw[warc] (-3, -1) -- (3, -1);
            \draw[warc] (-3, -1) to[out=west, in=west] (-3, -6);
            \draw[warc] (-3, -6) -- (3, -6);
            \draw[warc] (3, -6) to[out=east, in=east] (3, -1);
            \node[vJW, scale=0.1, label={[xshift=0pt, yshift=-1pt]above:\scriptsize $n$}] at (0,2) {};
            \node at (0, -5) [below] {\scriptsize $k$};
        \end{scope}
    }
}\mkern-12mu} = (-1)^k \frac{[n+1]}{[n-k+1]}\kbrack{\mkern-12mu\mbox{
    \tikz[baseline=-.6ex, scale=.075]{
        \setsquarecoords{12}
         
        \begin{scope}
            \draw[warc] (-12, 0) -- (-4, 0);
            \draw[warc] (-4, 0) -- (4, 0);
            \node[vJW, scale=0.1, label={[xshift=0pt, yshift=-1pt]above:\scriptsize $n-k$}] at (-4,0) {};
        \end{scope}
    }
}\mkern-12mu}$
    \end{enumerate}
    For properties (2) and (3), the property holds for the oppositely oriented twist by substituting $q \ra q^{-1}$.
\section{Colored Skein Relations}
    Now that Jones-Wenzl idempotents have been defined, we are able to color any tangle diagram, and use the above properties to reduce. Our main goal remains to compute colored Jones polynomials of links whose components carry different colors. Yuasa \cite{yuasa} gives the following twist formulas, which expand the Kauffman bracket of a half or full twist of two $n$-colored strands in a basis of untwisted diagrams. In the following section we generalize the full twist formulas to two strands colored $m$ and $n$.\
    
    \hspace{-0.25in}\textbf{Proposition 5.1 $\llbracket$Colored Kauffman Bracket Relation$\rrbracket$:} \cite[Prop.~3.1]{yuasa}
    For $n \in \zz^+$, and two colored $n-$strands crossing each other in a half-twist, we can resolve them in the following way:$$
\Biggl \langle \mkern-8mu \mbox{
    \tikz[baseline=-.6ex, scale=.08]{
        \setsquarecoords{12}
         
        \begin{scope}
             
            \draw[warc] (-8,6) to[out=east, in=west] (8,-6);
            \draw[warc] (-8,-6) to[out=east, in=west] (8,6);

            \draw[warc] (-12, 6) -- (-8, 6) node[midway, above, black, yshift=0pt, xshift=-2pt]{\scriptsize $n$};
            
            \draw[warc] (-12, -6) -- (-8, -6) node[midway, below, black, yshift=0pt, xshift=-2pt]{\scriptsize $n$};
            
            \draw[warc] (12, 6) -- (8, 6) node[midway, above, black, yshift=0pt, xshift=2pt]{\scriptsize $n$};
            
            \draw[warc] (12, -6) -- (8, -6) node[midway, below, black, yshift=0pt, xshift=2pt]{\scriptsize $n$};

            \node[vJW, scale=0.06] at (-8,6) {};
            \node[vJW, scale=0.06] at (-8,-6) {};
            \node[vJW, scale=0.06] at (8,6) {};
            \node[vJW, scale=0.06] at (8,-6) {};
        \end{scope}
    }
}\mkern-8mu \Biggr \rangle = \sum^n_{k=0} q^{\frac14 (-n^2+2k^2)} \binom{n}{k}_q 
\Biggl \langle \mkern-15mu \mbox{
    \tikz[baseline=-.6ex, scale=.08]{
        \setsquarecoords{12}
         
        \begin{scope}
             
            \draw[warc] (-8,-8) -- (8,-8) node[midway, below, black, yshift=0pt] {\scriptsize $k$};
            \draw[warc] (-8,8) -- (8,8) node[midway, above, black, yshift=0pt] {\scriptsize $k$};

            \draw[warc] (-8, 4) to[out=east, in =east, looseness=1.5, edge node = {node[ left, midway, black, font=\scriptsize] {$n-k$}}] (-8, -4);
             
            \draw[warc] (8, 4) to[out=west, in =west, looseness=1.5, edge node = {node[ right, midway, black, font=\scriptsize] {$n-k$}}] (8, -4);

            \draw[warc] (-12, 6) -- (-8, 6) node[midway, above, black, yshift=0pt, xshift=-1pt]{\scriptsize $n$};
            
            \draw[warc] (-12, -6) -- (-8, -6) node[midway, below, black, yshift=0pt, xshift=-1pt]{\scriptsize $n$};
            
            \draw[warc] (12, 6) -- (8, 6) node[midway, above, black, yshift=0pt, xshift=1pt]{\scriptsize $n$};
            
            \draw[warc] (12, -6) -- (8, -6) node[midway, below, black, yshift=0pt, xshift=1pt]{\scriptsize $n$};

            \node[vJW, scale=0.06] at (-8,6) {};
            \node[vJW, scale=0.06] at (-8,-6) {};
            \node[vJW, scale=0.06] at (8,6) {};
            \node[vJW, scale=0.06] at (8,-6) {};
        \end{scope}
    }
}\mkern-15mu \Biggr \rangle$$

    \hspace{-0.25in}\textbf{Proposition 5.2 $\llbracket$Full Twist Formula$\rrbracket$:} \cite[Prop.~3.3]{yuasa}
    For $n \in \zz^+$, and two colored $n-$strands crossing each other in a full twist, we can resolve them in the following way:
    $$
\Biggl \langle \mkern-8mu \mbox{
    \tikz[baseline=-.6ex, scale=.08]{
        \setsquarecoords{12}
         
        \begin{scope}
             
            \draw[warc] (-8,-6) to[out=east, in=west] (0,6);
            \draw[warc] (-8,6) to[out=east, in=west] (0,-6);
            \draw[warc] (0,-6) to[out=east, in=west] (8,6);
            \draw[warc] (0,6) to[out=east, in=west] (8,-6);

            \draw[warc] (-12, 6) -- (-8, 6) node[midway, above, black, yshift=0pt, xshift=-2pt]{\scriptsize $n$};
            
            \draw[warc] (-12, -6) -- (-8, -6) node[midway, below, black, yshift=0pt, xshift=-2pt]{\scriptsize $n$};
            
            \draw[warc] (12, 6) -- (8, 6) node[midway, above, black, yshift=0pt, xshift=2pt]{\scriptsize $n$};
            
            \draw[warc] (12, -6) -- (8, -6) node[midway, below, black, yshift=0pt, xshift=2pt]{\scriptsize $n$};

            \node[vJW, scale=0.06] at (-8,6) {};
            \node[vJW, scale=0.06] at (-8,-6) {};
            \node[vJW, scale=0.06] at (8,6) {};
            \node[vJW, scale=0.06] at (8,-6) {};
        \end{scope}
    }
}\mkern-8mu \Biggr \rangle  = \sum^n_{k=0}(-1)^{n-k}q^{{\frac12}(2k^2 - n^2 +k - n)}\frac{(q)_n}{(q)_k} \binom{n}{k}_q
\Biggl \langle \mkern-8mu
\mbox{
    \tikz[baseline=-.6ex, scale=.08]{
        \setsquarecoords{12}
         
        \begin{scope}
             
            \draw[warc] (-8,8) -- (8,8) node[midway, above, black, yshift=0pt] {\scriptsize $k$};
            \draw[warc] (-8,-8) -- (8,-8) node[midway, below, black, yshift=0pt] {\scriptsize $k$};

            \draw[warc] (-8, 4) to[out=east, in =east, looseness=1.5, edge node = {node[ left, midway, black, font=\scriptsize] {$n-k$}}] (-8, -4);
             
            \draw[warc] (8, 4) to[out=west, in =west, looseness=1.5, edge node = {node[ right, midway, black, font=\scriptsize] {$n-k$}}] (8, -4);

            \draw[warc] (-12, 6) -- (-8, 6) node[midway, above, black, yshift=0pt, xshift=-1pt]{\scriptsize $n$};
            
            \draw[warc] (-12, -6) -- (-8, -6) node[midway, below, black, yshift=0pt, xshift=-1pt]{\scriptsize $n$};
            
            \draw[warc] (12, 6) -- (8, 6) node[midway, above, black, yshift=0pt, xshift=1pt]{\scriptsize $n$};
            
            \draw[warc] (12, -6) -- (8, -6) node[midway, below, black, yshift=0pt, xshift=1pt]{\scriptsize $n$};

            \node[vJW, scale=0.06] at (-8,6) {};
            \node[vJW, scale=0.06] at (-8,-6) {};
            \node[vJW, scale=0.06] at (8,6) {};
            \node[vJW, scale=0.06] at (8,-6) {};
        \end{scope}
    }
} \mkern-8mu \Biggr \rangle$$
Both of these are used extensively in the reduction of colored tangle diagrams, and a proof can be given using the above properties of Jones-Wenzl idempotents, as they are in \cite{yuasa}. In fact, Yuasa also goes on to give a formula for any number of twists as follows:\

\hspace{-0.25in}\textbf{Proposition 5.3 $\llbracket$$p$-Half Twist Formula$\rrbracket$:} \cite[Prop.~3.5]{yuasa}
Given $n \in \zz^+$ and two of these strands crossing in $p-$half twists, we can resolve them in the following way, where $k_{i+1}' = k_{i} - k_{i+1}$ for $i = 0, 1, \cdots, p-1$.
$$\left\langle\mbox{
    \tikz[baseline=-.6ex, scale=.09]{
        \setsquarecoords{12}
         
        \begin{scope}
             
            \draw[warc] (-8,-6) to[out=east, in=west] (0,6);
            \draw[warc] (-8,6) to[out=east, in=west] (0,-6);

            \draw[warc] (-12, 6) -- (-8, 6) node[midway, above, black, yshift=0pt, xshift=-2pt]{\scriptsize $n$};
            
            \draw[warc] (-12, -6) -- (-8, -6) node[midway, below, black, yshift=0pt, xshift=-2pt]{\scriptsize $n$};

            \node[vJW, scale=0.06] at (-8,6) {};
            \node[vJW, scale=0.06] at (-8,-6) {};
            \node at (10,0) [align=center] {\scriptsize $...$};
             \node at (10,0) [above, align=center] {\scriptsize $p$ half\\[-8pt] \scriptsize twists};
        \end{scope}

        \begin{scope}[shift={(20,0)}]
            \draw[warc] (0,-6) to[out=east, in=west] (8,6);
            \draw[warc] (0,6) to[out=east, in=west] (8,-6);

            \draw[warc] (8, 6) -- (12, 6) node[midway, above, black, yshift=0pt, xshift=2pt]{\scriptsize $n$};
            
            \draw[warc] (8, -6) -- (12, -6) node[midway, below, black, yshift=0pt, xshift=2pt]{\scriptsize $n$};

            \node[vJW, scale=0.06] at (8,6) {};
            \node[vJW, scale=0.06] at (8,-6) {};
            
        \end{scope}
    }
}\right\rangle = (-1)^{pn}q^{-p(n^2 + 2n)/4} $$

$$\times \ \ \ \sum_{\substack{0 \leq k_p \leq \dots   \\  \leq k_1 \leq k_0 = n}}(-1)^{n-k_p}q^{(n-k_p)/2}(-1)^{\sum_{i=1}^p k_i}q^{\frac12 \sum_{i=1}^p (k_i^2 + k_i)} \binom{n}{k_1', k_2', \cdots, k_p', k_p}_q \Biggl \langle \mkern-8mu
\mbox{
    \tikz[baseline=-.6ex, scale=.08]{
        \setsquarecoords{12}
         
        \begin{scope}
             
            \draw[warc] (-8,8) -- (8,8) node[midway, above, black, yshift=0pt] {\scriptsize $k$};
            \draw[warc] (-8,-8) -- (8,-8) node[midway, below, black, yshift=0pt] {\scriptsize $k$};

            \draw[warc] (-8, 4) to[out=east, in =east, looseness=1.5, edge node = {node[ left, midway, black, font=\scriptsize] {$n-k$}}] (-8, -4);
             
            \draw[warc] (8, 4) to[out=west, in =west, looseness=1.5, edge node = {node[ right, midway, black, font=\scriptsize] {$n-k$}}] (8, -4);

            \draw[warc] (-12, 6) -- (-8, 6) node[midway, above, black, yshift=0pt, xshift=-1pt]{\scriptsize $n$};
            
            \draw[warc] (-12, -6) -- (-8, -6) node[midway, below, black, yshift=0pt, xshift=-1pt]{\scriptsize $n$};
            
            \draw[warc] (12, 6) -- (8, 6) node[midway, above, black, yshift=0pt, xshift=1pt]{\scriptsize $n$};
            
            \draw[warc] (12, -6) -- (8, -6) node[midway, below, black, yshift=0pt, xshift=1pt]{\scriptsize $n$};

            \node[vJW, scale=0.06] at (-8,6) {};
            \node[vJW, scale=0.06] at (-8,-6) {};
            \node[vJW, scale=0.06] at (8,6) {};
            \node[vJW, scale=0.06] at (8,-6) {};
        \end{scope}
    }
} \mkern-8mu \Biggr \rangle$$

\hspace{-0.25in}\textbf{Proposition 5.4 $\llbracket$$p$-Full Twist Formula$\rrbracket$:} \cite[Prop.~3.6]{yuasa}
Similarly for $n \in \zz^+$ and two of these colored strands crossing in $p-$full twists, we can resolve them in the following way where $k_{i+1}' = k_i - k_{i+1}$
$$\left\langle\mbox{
    \tikz[baseline=-.6ex, scale=.09]{
        \setsquarecoords{12}
         
        \begin{scope}
             
            \draw[warc] (-8,-6) to[out=east, in=west] (0,6);
            \draw[warc] (-8,6) to[out=east, in=west] (0,-6);
            \draw[warc] (0,-6) to[out=east, in=west] (8,6);
            \draw[warc] (0,6) to[out=east, in=west] (8,-6);

            \draw[warc] (-12, 6) -- (-8, 6) node[midway, above, black, yshift=0pt, xshift=-2pt]{\scriptsize $n$};
            
            \draw[warc] (-12, -6) -- (-8, -6) node[midway, below, black, yshift=0pt, xshift=-2pt]{\scriptsize $n$};

            \node[vJW, scale=0.06] at (-8,6) {};
            \node[vJW, scale=0.06] at (-8,-6) {};
            \node at (20,0) [align=center] {\scriptsize $...$};
             \node at (20,0) [above, align=center] {\scriptsize $p$ full\\[-8pt] \scriptsize twists};
        \end{scope}

        \begin{scope}[shift={(40,0)}]
            \draw[warc] (-8,-6) to[out=east, in=west] (0,6);
            \draw[warc] (-8,6) to[out=east, in=west] (0,-6);
            \draw[warc] (0,-6) to[out=east, in=west] (8,6);
            \draw[warc] (0,6) to[out=east, in=west] (8,-6);

            \draw[warc] (8, 6) -- (12, 6) node[midway, above, black, yshift=0pt, xshift=2pt]{\scriptsize $n$};
            
            \draw[warc] (8, -6) -- (12, -6) node[midway, below, black, yshift=0pt, xshift=2pt]{\scriptsize $n$};

            \node[vJW, scale=0.06] at (8,6) {};
            \node[vJW, scale=0.06] at (8,-6) {};
            
        \end{scope}
    }
}\right\rangle = q^{-p(n^2 + 2n)/2} \sum_{\substack{0 \leq k_p \leq \dots   \\  \leq k_1 \leq k_0 = n}}(-1)^{n-k_p}q^{(n-k_p)/2}$$

$$\times \ \ \ q^{\sum_{i=1}^p (k_i^2 + k_i)}\frac{(q)_n}{(q)_{k_p}} \binom{n}{k_1', k_2', \cdots, k_p', k_p}_q \Biggl \langle \mkern-8mu
\mbox{
    \tikz[baseline=-.6ex, scale=.08]{
        \setsquarecoords{12}
         
        \begin{scope}
             
            \draw[warc] (-8,8) -- (8,8) node[midway, above, black, yshift=0pt] {\scriptsize $k$};
            \draw[warc] (-8,-8) -- (8,-8) node[midway, below, black, yshift=0pt] {\scriptsize $k$};

            \draw[warc] (-8, 4) to[out=east, in =east, looseness=1.5, edge node = {node[ left, midway, black, font=\scriptsize] {$n-k$}}] (-8, -4);
             
            \draw[warc] (8, 4) to[out=west, in =west, looseness=1.5, edge node = {node[ right, midway, black, font=\scriptsize] {$n-k$}}] (8, -4);

            \draw[warc] (-12, 6) -- (-8, 6) node[midway, above, black, yshift=0pt, xshift=-1pt]{\scriptsize $n$};
            
            \draw[warc] (-12, -6) -- (-8, -6) node[midway, below, black, yshift=0pt, xshift=-1pt]{\scriptsize $n$};
            
            \draw[warc] (12, 6) -- (8, 6) node[midway, above, black, yshift=0pt, xshift=1pt]{\scriptsize $n$};
            
            \draw[warc] (12, -6) -- (8, -6) node[midway, below, black, yshift=0pt, xshift=1pt]{\scriptsize $n$};

            \node[vJW, scale=0.06] at (-8,6) {};
            \node[vJW, scale=0.06] at (-8,-6) {};
            \node[vJW, scale=0.06] at (8,6) {};
            \node[vJW, scale=0.06] at (8,-6) {};
        \end{scope}
    }
} \mkern-8mu \Biggr \rangle$$
\section{Twist Formulas for \texorpdfstring{$(m,n)-$}{(m,n)-}Colored Strands}
The result at the heart of this paper is a formula for $p$-full twists of two differently colored strands, labeled $m$ and $n$, which may be applied to link diagrams of differently colored components. For the purposes of this paper and its application to 2-bridge links in the next section, we focus only on the generalized colored full-twist formula. Setting $m=n$ in Proposition 6.1 and Theorem 6.3 recovers the full twist formulas of Propositions 5.2 and 5.4; the half twist formulas (Propositions 5.1 and 5.3) are not generalized here. Let us build up to the proof of the $p$-full twist formulas of $m,n$-colored strands.

\hspace{-0.25in}\textbf{Proposition 6.1: Main Result - $\llbracket$$(m,n)$-Colored Full Twist Formula$\rrbracket$:}
Given $n, m \in \zz^+$, and assuming $m \leq n$, we have the following formula for the $(m,n)-$colored full twist.
$$
\Biggl \langle \mkern-8mu \mbox{
    \tikz[baseline=-.6ex, scale=.08]{
        \setsquarecoords{12}
         
        \begin{scope}
             
            \draw[warc] (-8,-6) to[out=east, in=west] (0,6);
            \draw[warc] (-8,6) to[out=east, in=west] (0,-6);
            \draw[warc] (0,-6) to[out=east, in=west] (8,6);
            \draw[warc] (0,6) to[out=east, in=west] (8,-6);

            \draw[warc] (-12, 6) -- (-8, 6) node[midway, above, black, yshift=0pt, xshift=-2pt]{\scriptsize $n$};
            
            \draw[warc] (-12, -6) -- (-8, -6) node[midway, below, black, yshift=0pt, xshift=-2pt]{\scriptsize $m$};
            
            \draw[warc] (12, 6) -- (8, 6) node[midway, above, black, yshift=0pt, xshift=2pt]{\scriptsize $n$};
            
            \draw[warc] (12, -6) -- (8, -6) node[midway, below, black, yshift=0pt, xshift=2pt]{\scriptsize $m$};

            \node[vJW, scale=0.06] at (-8,6) {};
            \node[vJW, scale=0.06] at (-8,-6) {};
            \node[vJW, scale=0.06] at (8,6) {};
            \node[vJW, scale=0.06] at (8,-6) {};
        \end{scope}
    }
}\mkern-8mu \Biggr \rangle  = \sum^m_{k=0}(-1)^{m-k}q^{{\frac12}(2k(k+n-m) + k - m -nm)}\frac{(q)_n}{(q)_{k+n-m}} \binom{m}{k}_q
\Biggl \langle \mkern-8mu
\mbox{
    \tikz[baseline=-.6ex, scale=.08]{
        \setsquarecoords{12}
         
        \begin{scope}
             
            \draw[warc] (-8,8) -- (8,8) node[midway, above, black, yshift=0pt] {\tiny $k+n-m$};
            \draw[warc] (-8,-8) -- (8,-8) node[midway, below, black, yshift=0pt] {\scriptsize $k$};

            \draw[warc] (-8, 4) to[out=east, in =east, looseness=1.5, edge node = {node[ left, midway, black, font=\scriptsize] {$m-k$}}] (-8, -4);
             
            \draw[warc] (8, 4) to[out=west, in =west, looseness=1.5, edge node = {node[ right, midway, black, font=\scriptsize] {$m-k$}}] (8, -4);

            \draw[warc] (-12, 6) -- (-8, 6) node[midway, above, black, yshift=0pt, xshift=-1pt]{\scriptsize $n$};
            
            \draw[warc] (-12, -6) -- (-8, -6) node[midway, below, black, yshift=0pt, xshift=-2pt]{\scriptsize $m$};
            
            \draw[warc] (12, 6) -- (8, 6) node[midway, above, black, yshift=0pt, xshift=1pt]{\scriptsize $n$};
            
            \draw[warc] (12, -6) -- (8, -6) node[midway, below, black, yshift=0pt, xshift=2pt]{\scriptsize $m$};

            \node[vJW, scale=0.06] at (-8,6) {};
            \node[vJW, scale=0.06] at (-8,-6) {};
            \node[vJW, scale=0.06] at (8,6) {};
            \node[vJW, scale=0.06] at (8,-6) {};
        \end{scope}
    }
} \mkern-8mu \Biggr \rangle$$

\hspace{-0.25in}\textbf{Proof:} This follows the proof of \cite[Lemma~4.1]{masbaum}. Begin with the case $m =1$, and reduce using the Kauffman bracket to recover a recursion relation. First we begin by writing: 
$$
\Biggl \langle \mkern-8mu \mbox{
    \tikz[baseline=-.6ex, scale=.08]{
        \setsquarecoords{12}
         
        \begin{scope}
             
            \draw[warc] (-8,-6) to[out=east, in=west] (0,6);
            \draw[warc] (-8,6) to[out=east, in=west] (0,-6);
            \draw[warc] (0,-6) to[out=east, in=west] (8,6);
            \draw[warc] (0,6) to[out=east, in=west] (8,-6);

            \draw[warc] (-12, 6) -- (-8, 6) node[midway, above, black, yshift=0pt, xshift=-2pt]{\scriptsize $n$};
            
            \draw[warc] (-12, -6) -- (-8, -6) node[midway, below, black, yshift=0pt, xshift=-2pt]{\scriptsize $m$};
            
            \draw[warc] (12, 6) -- (8, 6) node[midway, above, black, yshift=0pt, xshift=2pt]{\scriptsize $n$};
            
            \draw[warc] (12, -6) -- (8, -6) node[midway, below, black, yshift=0pt, xshift=2pt]{\scriptsize $m$};

            \node[vJW, scale=0.06] at (-8,6) {};
            \node[vJW, scale=0.06] at (-8,-6) {};
            \node[vJW, scale=0.06] at (8,6) {};
            \node[vJW, scale=0.06] at (8,-6) {};
        \end{scope}
    }
}\mkern-8mu \Biggr \rangle  = \sum^m_{k=0} C_{n,m,k} \ 
\Biggl \langle \mkern-8mu
\mbox{
    \tikz[baseline=-.6ex, scale=.08]{
        \setsquarecoords{12}
         
        \begin{scope}
             
            \draw[warc] (-8,8) -- (8,8) node[midway, above, black, yshift=0pt] {\scriptsize $n-k$};
            \draw[warc] (-8,-8) -- (8,-8) node[midway, below, black, yshift=0pt] {\scriptsize $m-k$};

            \draw[warc] (-8, 4) to[out=east, in =east, looseness=1.5, edge node = {node[ left, midway, black, font=\scriptsize] {$k$}}] (-8, -4);
             
            \draw[warc] (8, 4) to[out=west, in =west, looseness=1.5, edge node = {node[ right, midway, black, font=\scriptsize] {$k$}}] (8, -4);

            \draw[warc] (-12, 6) -- (-8, 6) node[midway, above, black, yshift=0pt, xshift=-1pt]{\scriptsize $n$};
            
            \draw[warc] (-12, -6) -- (-8, -6) node[midway, below, black, yshift=0pt, xshift=-2pt]{\scriptsize $m$};
            
            \draw[warc] (12, 6) -- (8, 6) node[midway, above, black, yshift=0pt, xshift=1pt]{\scriptsize $n$};
            
            \draw[warc] (12, -6) -- (8, -6) node[midway, below, black, yshift=0pt, xshift=2pt]{\scriptsize $m$};

            \node[vJW, scale=0.06] at (-8,6) {};
            \node[vJW, scale=0.06] at (-8,-6) {};
            \node[vJW, scale=0.06] at (8,6) {};
            \node[vJW, scale=0.06] at (8,-6) {};
        \end{scope}
    }
} \mkern-8mu \Biggr \rangle$$
$$
\Biggl \langle \mkern-8mu \mbox{
    \tikz[baseline=-.6ex, scale=.08]{
        \setsquarecoords{12}
         
        \begin{scope}
             
            \draw[warc] (-8,-6) to[out=east, in=west] (0,6);
            \draw[owarc] (-8,6) to[out=east, in=west] (0,-6);
            \draw[owarc] (0,-6) to[out=east, in=west] (8,6);
            \draw[warc] (0,6) to[out=east, in=west] (8,-6);

            \draw[warc] (-12, 6) -- (-8, 6) node[midway, above, black, yshift=0pt, xshift=-2pt]{\scriptsize $n$};
            
            \draw[warc] (-12, -6) -- (-8, -6) node[midway, below, black, yshift=0pt, xshift=-2pt]{\scriptsize $1$};
            
            \draw[warc] (12, 6) -- (8, 6) node[midway, above, black, yshift=0pt, xshift=2pt]{\scriptsize $n$};
            
            \draw[warc] (12, -6) -- (8, -6) node[midway, below, black, yshift=0pt, xshift=2pt]{\scriptsize $1$};

            \node at (0,-12) [above, align=center] {\scriptsize $n-1$};

            \node at (0,-6) [above, align=center] {\scriptsize $1$};

            \node[vJW, scale=0.06] at (-8,6) {};
            \node[vJW, scale=0.06] at (8,6) {};
        \end{scope}
    }
}\mkern-8mu \Biggr \rangle = q^{1/2} \Biggl \langle \mkern-8mu \mbox{
    \tikz[baseline=-.6ex, scale=.08]{
        \setsquarecoords{12}
         
        \begin{scope}

            \draw[warc] (-8, 8) -- (8, 8);
            \draw[warc] (-8,-6) to[out=east, in=west] (0,6);
            \draw[warc] (-8,6) to[out=east, in=west] (0,-6);
            \draw[warc] (0,-6) to[out=east, in=west] (8,6);
            \draw[warc] (0,6) to[out=east, in=west] (8,-6);

            \draw[warc] (-12, 6) -- (-8, 6) node[midway, above, black, yshift=0pt, xshift=-2pt]{\scriptsize $n$};
            
            \draw[warc] (-12, -6) -- (-8, -6) node[midway, below, black, yshift=0pt, xshift=-2pt]{\scriptsize $1$};
            
            \draw[warc] (12, 6) -- (8, 6) node[midway, above, black, yshift=0pt, xshift=2pt]{\scriptsize $n$};
            
            \draw[warc] (12, -6) -- (8, -6) node[midway, below, black, yshift=0pt, xshift=2pt]{\scriptsize $1$};

            \node at (0,-12) [above, align=center] {\scriptsize $n-1$};

            \node at (0,0) [above, align=center] {\scriptsize $1$};

            \node at (0,7) [above, align=center] {\scriptsize $1$};

            \node[vJW, scale=0.06] at (-8,6) {};
            \node[vJW, scale=0.06] at (8,6) {};
        \end{scope}
    }
}\mkern-8mu \Biggr \rangle + q^{-(n-1)/2}(1-q^{-1}) \Biggl \langle \mkern-8mu \mbox{
    \tikz[baseline=-.6ex, scale=.08]{
        \setsquarecoords{12}
         
        \begin{scope}

            \draw[warc] (-8, 8) -- (8, 8);
            \draw[warc] (-8,-6) to[out=east, in=east] (-8,6);
            \draw[warc] (8,6) to[out=west, in=west] (8,-6);

            \draw[warc] (-12, 6) -- (-8, 6) node[midway, above, black, yshift=0pt, xshift=-2pt]{\scriptsize $n$};
            
            \draw[warc] (-12, -6) -- (-8, -6) node[midway, below, black, yshift=0pt, xshift=-2pt]{\scriptsize $1$};
            
            \draw[warc] (12, 6) -- (8, 6) node[midway, above, black, yshift=0pt, xshift=2pt]{\scriptsize $n$};
            
            \draw[warc] (12, -6) -- (8, -6) node[midway, below, black, yshift=0pt, xshift=2pt]{\scriptsize $1$};

            \node at (0,7) [above, align=center] {\scriptsize $n-1$};

            \node[vJW, scale=0.06] at (-8,6) {};
            \node[vJW, scale=0.06] at (8,6) {};
        \end{scope}
    }
}\mkern-8mu \Biggr \rangle$$
And it can be shown using induction on $m$, that:
$$\implies C_{n,m+1,k} = q^{\frac12(n-2k)}C_{n,m,k} + q^{\frac12(n - 2k+1)}(1 - q^{-n+k-1})C_{n,m, {k-1}}$$
for which $$C_{n,m,k} = q^{\frac12 m(n-k)}\bracketbinom{m}{k}_q \prod_{j=n-k+1}^n(1 - q^{-j})$$ is the unique solution. In the case $m=n$ this is \cite[Lemma~4.1]{masbaum}; the case $m \neq n$ is new.  This solution must be unique due to the 'binomial' like recurrence property above. That is, that $C_{n,0,0} = 1$ and $C_{n,m,k}=0$ for $k <0$ or $k>m$, so the product only needs to go up to $m$. We can rewrite $C_{n,m,k}$ to better fit our structure, using the following properties:
$$\bracketbinom{m}{k}_q = q^{-(m-k)(k)/2}\binom{m}{k}_q \quad (q')_n = (-1)^n q^{-n(n+1)/2} (q)_n$$
We will also make the substitution $k \ra m-k$ (with $q' = q^{-1}$), in order to match the form of Proposition 6.1.
$$\implies C_{n,m,k} = q^{m(k + n - m)/2} \bracketbinom{m}{m-k}_q \cdot \frac{(q')_n}{(q')_{k+n - m}}$$
$$ C_{n,m,k} = q^{m(k+n-m)/2}\cdot q^{-(m-k)k/2} \cdot\binom{m}{m-k}_q \cdot \frac{(-1)^nq^{-n(n+1)/2}}{(-1)^{k+n-m}q^{-(k+n-m)(k+n-m+1)/2}}\cdot \frac{(q)_n}{(q)_{k+n-m}}$$
$$\implies \Biggl \langle \mkern-8mu \mbox{
    \tikz[baseline=-.6ex, scale=.08]{
        \setsquarecoords{12}
         
        \begin{scope}
             
            \draw[warc] (-8,-6) to[out=east, in=west] (0,6);
            \draw[warc] (-8,6) to[out=east, in=west] (0,-6);
            \draw[warc] (0,-6) to[out=east, in=west] (8,6);
            \draw[warc] (0,6) to[out=east, in=west] (8,-6);

            \draw[warc] (-12, 6) -- (-8, 6) node[midway, above, black, yshift=0pt, xshift=-2pt]{\scriptsize $n$};
            
            \draw[warc] (-12, -6) -- (-8, -6) node[midway, below, black, yshift=0pt, xshift=-2pt]{\scriptsize $m$};
            
            \draw[warc] (12, 6) -- (8, 6) node[midway, above, black, yshift=0pt, xshift=2pt]{\scriptsize $n$};
            
            \draw[warc] (12, -6) -- (8, -6) node[midway, below, black, yshift=0pt, xshift=2pt]{\scriptsize $m$};

            \node[vJW, scale=0.06] at (-8,6) {};
            \node[vJW, scale=0.06] at (-8,-6) {};
            \node[vJW, scale=0.06] at (8,6) {};
            \node[vJW, scale=0.06] at (8,-6) {};
        \end{scope}
    }
}\mkern-8mu \Biggr \rangle  = \sum^m_{k=0}(-1)^{m-k}q^{{\frac12}(2k(k+n-m) + k - m -nm)}\frac{(q)_n}{(q)_{k+n-m}} \binom{m}{k}_q
\Biggl \langle \mkern-8mu
\mbox{
    \tikz[baseline=-.6ex, scale=.08]{
        \setsquarecoords{12}
         
        \begin{scope}
             
            \draw[warc] (-8,8) -- (8,8) node[midway, above, black, yshift=0pt] {\tiny $k+n-m$};
            \draw[warc] (-8,-8) -- (8,-8) node[midway, below, black, yshift=0pt] {\scriptsize $k$};

            \draw[warc] (-8, 4) to[out=east, in =east, looseness=1.5, edge node = {node[ left, midway, black, font=\scriptsize] {$m-k$}}] (-8, -4);
             
            \draw[warc] (8, 4) to[out=west, in =west, looseness=1.5, edge node = {node[ right, midway, black, font=\scriptsize] {$m-k$}}] (8, -4);

            \draw[warc] (-12, 6) -- (-8, 6) node[midway, above, black, yshift=0pt, xshift=-1pt]{\scriptsize $n$};
            
            \draw[warc] (-12, -6) -- (-8, -6) node[midway, below, black, yshift=0pt, xshift=-2pt]{\scriptsize $m$};
            
            \draw[warc] (12, 6) -- (8, 6) node[midway, above, black, yshift=0pt, xshift=1pt]{\scriptsize $n$};
            
            \draw[warc] (12, -6) -- (8, -6) node[midway, below, black, yshift=0pt, xshift=2pt]{\scriptsize $m$};

            \node[vJW, scale=0.06] at (-8,6) {};
            \node[vJW, scale=0.06] at (-8,-6) {};
            \node[vJW, scale=0.06] at (8,6) {};
            \node[vJW, scale=0.06] at (8,-6) {};
        \end{scope}
    }
} \mkern-8mu \Biggr \rangle$$

\hfill$\square$\\
\textbf{Proposition 6.2 $\llbracket$Sliding Formula$\rrbracket$:} For the colored full twist, this formula allows us to slide over components through full twists, and ends up being instrumental in the proof of the $p$ $(m,n)$-colored twist formulas.

$$\Biggl \langle \mkern-8mu \mbox{
    \tikz[baseline=-.6ex, scale=.08]{
        \setsquarecoords{12}
         
        \begin{scope}
             
            \draw[warc] (-8,-6) to[out=east, in=west] (0,6);
            \draw[warc] (-8,6) to[out=east, in=west] (0,-6);
            \draw[warc] (0,-6) to[out=east, in=west] (8,6);
            \draw[warc] (0,6) to[out=east, in=west] (8,-6);

            \draw[warc] (-12, 6) -- (-8, 6);
            
            \draw[warc] (-12, -6) -- (-8, -6);
            
            \draw[warc] (28, 6) -- (8, 6);
            
            \draw[warc] (28, -6) -- (8, -6);

            \draw[warc] (8, 5) to[out=east, in=east, looseness=1.25] (8, -5);

            \node[vJW, scale=0.06] at (-8,6) {};
            \node[vJW, scale=0.06] at (-8,-6) {};
            \node[vJW, scale=0.06] at (8,6) {};
            \node[vJW, scale=0.06] at (8,-6) {};

            \node at (-10.5,6) [above, align=center] {\scriptsize $n$};
            \node at (-11,-6) [below, align=center] {\scriptsize $m$};
            \node at (11,0) [right, align=center] {\scriptsize $k$};
            \node at (20,-5.5) [below, align=center] {\scriptsize $m-k$};
            \node at (20,5.5) [above, align=center] {\scriptsize $n-k$};
        \end{scope}
    }
}\mkern-8mu \Biggr \rangle = \frac{\mu_{m-k}\mu_{n-k}}{\mu_m \mu_n} \kbrack{\mkern-12mu\mbox{
    \tikz[baseline=-.6ex, scale=.075]{
        \setsquarecoords{12}
         
        \begin{scope}
             
            \draw[warc] (-8,-6) to[out=east, in=west] (0,6);
            \draw[warc] (-8,6) to[out=east, in=west] (0,-6);
            \draw[warc] (0,-6) to[out=east, in=west] (8,6);
            \draw[warc] (0,6) to[out=east, in=west] (8,-6);

            \draw[warc] (-10, 4) -- (-8, 4);
            \draw[warc] (-8, 4) to[out=east, in=east] (-8, -4);
            \draw[warc] (-10, -4) -- (-8, -4);
            \node at (-9,0) [align=center] {\scriptsize $k$};
             
            \draw[warc] (-14, 6) -- (-8, 6);
            
            \draw[warc] (-14, -6) -- (-8, -6);

            \node[vJW, scale=0.06, label={[xshift=2pt, yshift=6pt]left:\scriptsize $n$}] at (-10,6) {};
            
            \node[vJW, scale=0.06, label={[xshift=2pt, yshift=-6pt]left:\scriptsize $m$}] at (-10,-6) {};

            \node at (0.5,-5) [below, align=center] {\scriptsize $n-k$};
            \node at (0.5,5) [above, align=center] {\scriptsize $m-k$};

            \draw[warc] (8, 6) -- (14, 6);
            \draw[warc] (8, -6) -- (14, -6);

            \node[vJW, scale=0.06, label={[xshift=-2pt, yshift=-6pt]right:\scriptsize $m$}] at (10,-6) {};
            \node[vJW, scale=0.06, label={[xshift=-2pt, yshift=6pt]right:\scriptsize $n$}] at (10,6) {};
        \end{scope}     
    }
}\mkern-12mu}$$

\hspace{-0.25in}\textbf{Proof:} We are able to slide the rightmost JW idempotents to the right side, and reduce using the properties in Proposition 4.5.

$$
    \Biggl \langle \mkern-8mu \mbox{
    \tikz[baseline=-.6ex, scale=.08]{
        \setsquarecoords{12}
         
        \begin{scope}
             
            \draw[warc] (-8,-6) to[out=east, in=west] (0,6);
            \draw[warc] (-8,6) to[out=east, in=west] (0,-6);
            \draw[warc] (0,-6) to[out=east, in=west] (8,6);
            \draw[warc] (0,6) to[out=east, in=west] (8,-6);

            \draw[warc] (-12, 6) -- (-8, 6);
            
            \draw[warc] (-12, -6) -- (-8, -6);
            
            \draw[warc] (28, 6) -- (8, 6);
            
            \draw[warc] (28, -6) -- (8, -6);

            \draw[warc] (8, 5) to[out=east, in=east, looseness=1.25] (8, -5);

            \node[vJW, scale=0.06] at (-8,6) {};
            \node[vJW, scale=0.06] at (-8,-6) {};
            \node[vJW, scale=0.06] at (8,6) {};
            \node[vJW, scale=0.06] at (8,-6) {};

            \node at (-10.5,6) [above, align=center] {\scriptsize $n$};
            \node at (-11,-6) [below, align=center] {\scriptsize $m$};
            \node at (11,0) [right, align=center] {\scriptsize $k$};
            \node at (20,-5.5) [below, align=center] {\scriptsize $m-k$};
            \node at (20,5.5) [above, align=center] {\scriptsize $n-k$};
        \end{scope}
    }
}\mkern-8mu \Biggr \rangle = \Biggl \langle \mkern-8mu \mbox{
    \tikz[baseline=-.6ex, scale=.08]{
        \setsquarecoords{12}
         
        \begin{scope}
             
            \draw[owarc] (-12, -5.5) to[out=east, in=west] (-2, 5.5);
            \draw[owarc] (-12, 5.5) to[out=east, in=west] (-2, -5.5);
            \draw[owarc] (-2, -5.5) to[out=east, in=west] (8, 5.5);
            \draw[owarc] (-2, 5.5) to[out=east, in=west] (8, -5.5);

            \draw[warc] (-16, 6) -- (-12, 6);
            
            \draw[warc] (-16, -6) -- (-12, -6);
            
            \draw[warc] (18, 6) -- (8, 6);
            
            \draw[warc] (18, -6) -- (8, -6);

            \draw[warc] (8, 5) to[out=east, in=east, looseness=1.25] (8, -5);

            \node[vJW, scale=0.06] at (-12,6) {};
            \node[vJW, scale=0.06] at (-12,-6) {};

            \node at (-14.5,6) [above, align=center] {\scriptsize $n$};
            \node at (-15,-6) [below, align=center] {\scriptsize $m$};
            \node at (11,0) [right, align=center] {\scriptsize $k$};
            \node at (14,-5.5) [below, align=center] {\scriptsize $m-k$};
            \node at (14,5.5) [above, align=center] {\scriptsize $n-k$};
        \end{scope}
    }
}\mkern-8mu \Biggr \rangle = \mu^{-1}_k q^{k(n-k)/4}q^{k(m-k)/4}\Biggl \langle \mkern-8mu \mbox{
    \tikz[baseline=-.6ex, scale=.08]{
        \setsquarecoords{12}
         
        \begin{scope}
             
            \draw[owarc] (-12, -6) to[out=east, in=west] (-2, 5.5);
            \draw[owarc] (-12, 6) to[out=east, in=west] (-2, -5.5);
            \draw[warc] (-2, -5) to[out=east, in=east] (-2, 5);

            \draw[warc] (-2, -6) to[out=east, in=west] (8, 6);

            \draw[warc] (-2, 6) to[out=east, in=west] (8, -6);

            \draw[warc] (-16, 6) -- (-12, 6);
            
            \draw[warc] (-16, -6) -- (-12, -6);
            
            \draw[warc] (18, 6) -- (8, 6);
            
            \draw[warc] (18, -6) -- (8, -6);

            \node[vJW, scale=0.06] at (-12,6) {};
            \node[vJW, scale=0.06] at (-12,-6) {};

            \node at (-14.5,6) [above, align=center] {\scriptsize $n$};
            \node at (-15,-6) [below, align=center] {\scriptsize $m$};
            \node at (-4,0) [right, align=center] {\scriptsize $k$};
            \node at (14,-5.5) [below, align=center] {\scriptsize $m-k$};
            \node at (14,5.5) [above, align=center] {\scriptsize $n-k$};
        \end{scope}
    }
}\mkern-8mu \Biggr \rangle $$
\\
$$= \mu_k^{-2} q^{-k(n-k)/2}q^{-k(m-k)/2}\kbrack{\mkern-12mu\mbox{
    \tikz[baseline=-.6ex, scale=.075]{
        \setsquarecoords{12}
         
        \begin{scope}
             
            \draw[warc] (-8,-6) to[out=east, in=west] (0,6);
            \draw[warc] (-8,6) to[out=east, in=west] (0,-6);
            \draw[warc] (0,-6) to[out=east, in=west] (8,6);
            \draw[warc] (0,6) to[out=east, in=west] (8,-6);

            \draw[warc] (-10, 4) -- (-8, 4);
            \draw[warc] (-8, 4) to[out=east, in=east] (-8, -4);
            \draw[warc] (-10, -4) -- (-8, -4);
            \node at (-9,0) [align=center] {\scriptsize $k$};
             
            \draw[warc] (-14, 6) -- (-8, 6);
            
            \draw[warc] (-14, -6) -- (-8, -6);

            \node[vJW, scale=0.06, label={[xshift=2pt, yshift=6pt]left:\scriptsize $n$}] at (-10,6) {};
            
            \node[vJW, scale=0.06, label={[xshift=2pt, yshift=-6pt]left:\scriptsize $m$}] at (-10,-6) {};

            \node at (0.5,-5) [below, align=center] {\scriptsize $n-k$};
            \node at (0.5,5) [above, align=center] {\scriptsize $m-k$};

            \draw[warc] (8, 6) -- (14, 6);
            \draw[warc] (8, -6) -- (14, -6);

            \node[vJW, scale=0.06, label={[xshift=-2pt, yshift=-6pt]right:\scriptsize $m$}] at (10,-6) {};
            \node[vJW, scale=0.06, label={[xshift=-2pt, yshift=6pt]right:\scriptsize $n$}] at (10,6) {};
        \end{scope}     
    }
}\mkern-12mu} = \frac{\mu_{m-k}\mu_{n-k}}{\mu_m \mu_n}\kbrack{\mkern-12mu\mbox{
    \tikz[baseline=-.6ex, scale=.075]{
        \setsquarecoords{12}
         
        \begin{scope}
             
            \draw[warc] (-8,-6) to[out=east, in=west] (0,6);
            \draw[warc] (-8,6) to[out=east, in=west] (0,-6);
            \draw[warc] (0,-6) to[out=east, in=west] (8,6);
            \draw[warc] (0,6) to[out=east, in=west] (8,-6);

            \draw[warc] (-10, 4) -- (-8, 4);
            \draw[warc] (-8, 4) to[out=east, in=east] (-8, -4);
            \draw[warc] (-10, -4) -- (-8, -4);
            \node at (-9,0) [align=center] {\scriptsize $k$};
             
            \draw[warc] (-14, 6) -- (-8, 6);
            
            \draw[warc] (-14, -6) -- (-8, -6);

            \node[vJW, scale=0.06, label={[xshift=2pt, yshift=6pt]left:\scriptsize $n$}] at (-10,6) {};
            
            \node[vJW, scale=0.06, label={[xshift=2pt, yshift=-6pt]left:\scriptsize $m$}] at (-10,-6) {};

            \node at (0.5,-5) [below, align=center] {\scriptsize $n-k$};
            \node at (0.5,5) [above, align=center] {\scriptsize $m-k$};

            \draw[warc] (8, 6) -- (14, 6);
            \draw[warc] (8, -6) -- (14, -6);

            \node[vJW, scale=0.06, label={[xshift=-2pt, yshift=-6pt]right:\scriptsize $m$}] at (10,-6) {};
            \node[vJW, scale=0.06, label={[xshift=-2pt, yshift=6pt]right:\scriptsize $n$}] at (10,6) {};
        \end{scope}     
    }
}\mkern-12mu}$$
Note that $q^{-k(n-k)/2} = \mu_k \mu_{n-k} \mu_n^{-1}$ and similarly, $q^{-k(m-k)/2} = \mu_k \mu_{m-k} \mu_m^{-1}$ is where the above substitution came from. \hfill$\square$\

\hspace{-0.25in}\textbf{Corollary:}
$$\Biggl \langle \mkern-8mu \mbox{
    \tikz[baseline=-.6ex, scale=.08]{
        \setsquarecoords{12}
         
        \begin{scope}
             
            \draw[warc] (-8,-6) to[out=east, in=west] (0,6);
            \draw[warc] (-8,6) to[out=east, in=west] (0,-6);

            \draw[warc] (-12, 6) -- (-8, 6);
            
            \draw[warc] (-12, -6) -- (-8, -6);
            
            \draw[warc] (20, 6) -- (0, 6);
            
            \draw[warc] (20, -6) -- (0, -6);

            \draw[warc] (0, 5) to[out=east, in=east, looseness=1.25] (0, -5);

            \node[vJW, scale=0.06] at (-8,6) {};
            \node[vJW, scale=0.06] at (-8,-6) {};
            \node[vJW, scale=0.06] at (0,6) {};
            \node[vJW, scale=0.06] at (0,-6) {};

            \node at (-10.5,6) [above, align=center] {\scriptsize $n$};
            \node at (-11,-6) [below, align=center] {\scriptsize $m$};
            \node at (3,0) [right, align=center] {\scriptsize $k$};
            \node at (10,-5.5) [below, align=center] {\scriptsize $m-k$};
            \node at (10,5.5) [above, align=center] {\scriptsize $n-k$};
        \end{scope}
    }
}\mkern-8mu \Biggr \rangle = \left(\frac{\mu_{m-k}\mu_{n-k}}{\mu_m \mu_n}\right)^{1/2} \kbrack{\mkern-12mu\mbox{
    \tikz[baseline=-.6ex, scale=.075]{
        \setsquarecoords{12}
         
        \begin{scope}
             
            \draw[warc] (-8,-6) to[out=east, in=west] (0,6);
            \draw[warc] (-8,6) to[out=east, in=west] (0,-6);
            \draw[warc] (0,6) -- (10,6);
            \draw[warc] (0,-6) -- (10, -6);

            \draw[warc] (-10, 4) -- (-8, 4);
            \draw[warc] (-8, 4) to[out=east, in=east] (-8, -4);
            \draw[warc] (-10, -4) -- (-8, -4);
            \node at (-9,0) [align=center] {\scriptsize $k$};
             
            \draw[warc] (-14, 6) -- (-8, 6);
            
            \draw[warc] (-14, -6) -- (-8, -6);

            \node[vJW, scale=0.06, label={[xshift=2pt, yshift=6pt]left:\scriptsize $n$}] at (-10,6) {};
            
            \node[vJW, scale=0.06, label={[xshift=2pt, yshift=-6pt]left:\scriptsize $m$}] at (-10,-6) {};

            \draw[warc] (8, 6) -- (14, 6);
            \draw[warc] (8, -6) -- (14, -6);

            \node[vJW, scale=0.06, label={[xshift=-2pt, yshift=-6pt]right:\scriptsize $n-k$}] at (6,-6) {};
            \node[vJW, scale=0.06, label={[xshift=-2pt, yshift=6pt]right:\scriptsize $m-k$}] at (6,6) {};
        \end{scope}     
    }
}\mkern-12mu}$$
This comes from one of the intermediate steps from the proof above, and is the half-twist analogue of the sliding formula.\

\hspace{-0.25in}\textbf{Theorem 6.3 $\llbracket$$p$ Full Twists of $(m,n)-$Colored Strands.$\rrbracket$:}
For $m, n \in \zz^+$, and assuming $m \leq n$, and for $k_{i+1}' = k_i - k_{i+1}$, we have the following formula for $p$-full twists of (m,n)-colored strands:
    $$\kbrack{\mbox{
    \tikz[baseline=-.6ex, scale=.075]{
        \setsquarecoords{12}
         
        \begin{scope}
             
            \draw[warc] (-8,-6) to[out=east, in=west] (0,6);
            \draw[warc] (-8,6) to[out=east, in=west] (0,-6);
            \draw[warc] (0,-6) to[out=east, in=west] (8,6);
            \draw[warc] (0,6) to[out=east, in=west] (8,-6);

            \draw[warc] (-12, 6) -- (-8, 6) node[midway, above, black, yshift=0pt, xshift=-2pt]{\scriptsize $n$};
            
            \draw[warc] (-12, -6) -- (-8, -6) node[midway, below, black, yshift=0pt, xshift=-2pt]{\scriptsize $m$};

            \node[vJW, scale=0.06] at (-8,6) {};
            \node[vJW, scale=0.06] at (-8,-6) {};
            \node at (20,0) [align=center] {\scriptsize $...$};
             \node at (20,0) [above, align=center] {\scriptsize $p$ full\\[-8pt] \scriptsize twists};
        \end{scope}

        \begin{scope}[shift={(40,0)}]
            \draw[warc] (-8,-6) to[out=east, in=west] (0,6);
            \draw[warc] (-8,6) to[out=east, in=west] (0,-6);
            \draw[warc] (0,-6) to[out=east, in=west] (8,6);
            \draw[warc] (0,6) to[out=east, in=west] (8,-6);

            \draw[warc] (8, 6) -- (12, 6) node[midway, above, black, yshift=0pt, xshift=2pt]{\scriptsize $n$};
            
            \draw[warc] (8, -6) -- (12, -6) node[midway, below, black, yshift=0pt, xshift=2pt]{\scriptsize $m$};

            \node[vJW, scale=0.06] at (8,6) {};
            \node[vJW, scale=0.06] at (8,-6) {};
            
        \end{scope}
    }
}}
= q^{\frac{-p}{2}(mn+2m)}\sum_{\substack{0 \leq k_p \leq \cdots \\ \leq k_1 \leq k_0 = m}} (-1)^{m - k_p} q^{\frac{m-k_p}{2}}$$
$$\times q^{\sum_{i=1}^p[k_i^2 + k_i(n-m+1)]} \frac{(q)_n}{(q)_{k_p + n - m}} \binom{m}{k_1', \cdots, k_{p}', k_p}_q \Biggl \langle \mkern-8mu
\mbox{
    \tikz[baseline=-.6ex, scale=.08]{
        \setsquarecoords{12}
         
        \begin{scope}
             
            \draw[warc] (-8,8) -- (8,8) node[midway, above, black, yshift=0pt] {\tiny $k_p+n-m$};
            \draw[warc] (-8,-8) -- (8,-8) node[midway, below, black, yshift=0pt] {\scriptsize $k_p$};

            \draw[warc] (-8, 4) to[out=east, in =east, looseness=1.5, edge node = {node[ left, midway, black, font=\scriptsize] {$m-k_p$}}] (-8, -4);
             
            \draw[warc] (8, 4) to[out=west, in =west, looseness=1.5, edge node = {node[ right, midway, black, font=\scriptsize] {$m-k_p$}}] (8, -4);

            \draw[warc] (-12, 6) -- (-8, 6) node[midway, above, black, yshift=0pt, xshift=-1pt]{\scriptsize $n$};
            
            \draw[warc] (-12, -6) -- (-8, -6) node[midway, below, black, yshift=0pt, xshift=-2pt]{\scriptsize $m$};
            
            \draw[warc] (12, 6) -- (8, 6) node[midway, above, black, yshift=0pt, xshift=1pt]{\scriptsize $n$};
            
            \draw[warc] (12, -6) -- (8, -6) node[midway, below, black, yshift=0pt, xshift=2pt]{\scriptsize $m$};

            \node[vJW, scale=0.06] at (-8,6) {};
            \node[vJW, scale=0.06] at (-8,-6) {};
            \node[vJW, scale=0.06] at (8,6) {};
            \node[vJW, scale=0.06] at (8,-6) {};
        \end{scope}
    }
} \mkern-8mu \Biggr \rangle$$

\hspace{-0.25in}\textbf{Proof:}
In order to evaluate the $p$ full twist formula of $m$ and $n$ colored strands, let us define the following diagram then proceed inductively:
$$\kbrack{\sigma_{m,n}^{p-i}(k_i)} = 
\kbrack{\mkern-12mu\mbox{
    \tikz[baseline=-.6ex, scale=.08]{
        \setsquarecoords{12}
         
        \begin{scope}
             
            \draw[warc] (-8,-6) to[out=east, in=west] (0,6);
            \draw[warc] (-8,6) to[out=east, in=west] (0,-6);
            \draw[warc] (0,-6) to[out=east, in=west] (8,6);
            \draw[warc] (0,6) to[out=east, in=west] (8,-6);

            \draw[warc] (-10, 4) -- (-8, 4);
            \draw[warc] (-8, 4) to[out=east, in=east] (-8, -4);
            \draw[warc] (-10, -4) -- (-8, -4);
            \node at (-13,0) [align=center] {\tiny $m-k_i$};
             
            \draw[warc] (-14, 6) -- (-8, 6);
            
            \draw[warc] (-14, -6) -- (-8, -6);

            \node[vJW, scale=0.06, label={[xshift=2pt, yshift=6pt]left:\scriptsize $n$}] at (-10,6) {};
            
            \node[vJW, scale=0.06, label={[xshift=2pt, yshift=-6pt]left:\scriptsize $m$}] at (-10,-6) {};
            
            \node at (20,0) [align=center] {\scriptsize $...$};
            \node at (20,0) [above, align=center] {\scriptsize $p-i$ full\\[-8pt] \scriptsize twists};

            \node at (0.5,-5) [below, align=center] {\tiny $k_i + n-m$};
            \node at (0.5,5) [above, align=center] {\tiny $k_i$};
        \end{scope}

        \begin{scope}[shift={(40,0)}]
            \draw[warc] (-8,-6) to[out=east, in=west] (0,6);
            \draw[warc] (-8,6) to[out=east, in=west] (0,-6);
            \draw[warc] (0,-6) to[out=east, in=west] (8,6);
            \draw[warc] (0,6) to[out=east, in=west] (8,-6);

            \draw[warc] (8, 6) -- (14, 6);
            \draw[warc] (8, -6) -- (14, -6);

            \draw[warc] (10, 4) -- (8, 4);
            \draw[warc] (8, 4) to[out=west, in=west] (8, -4);
            \draw[warc] (10, -4) -- (8, -4);
            \node at (13,0) [align=center] {\tiny $m-k_i$};

            \node[vJW, scale=0.06, label={[xshift=-2pt, yshift=-6pt]right:\scriptsize $m$}] at (10,-6) {};
            \node[vJW, scale=0.06, label={[xshift=-2pt, yshift=6pt]right:\scriptsize $n$}] at (10,6) {};

            \node at (-0.5,-5) [below, align=center] {\tiny $k_i + n-m$};
            \node at (0.5,5) [above, align=center] {\tiny $k_i$};
            
        \end{scope}

    }
}\mkern-12mu}$$
Notice that the diagram we are trying to reduce in the above theorem is $\kbrack{\sigma^{p}_{m,n}(k_0)}$ where $k_0 = m$. Defining it this way allows us to build our result up inductively. Using the sliding formula of Proposition 6.2 as well as the single (m,n)-colored full twist reduction of Proposition 6.1, we are able to reduce a full twist at the rightmost end of the diagram, then 'slide' its remaining components completely to the right by applying the sliding formula $p-i-1$ times in the following way:
$$\kbrack{\sigma_{m,n}^{p-i}(k_i)} = \sum_{k_{i+1}=0}^{k_i} q^{\frac12(2k_{i+1}^2 + 2k_{i+1}(n-m) - k_i(k_i + n - m) - k_i + k_{i+1})} (-1)^{k_i - k_{i+1}} \frac{(q)_{k_i + n - m} }{(q)_{k_{i+1} + n - m}} \binom{k_i}{k_i - k_{i+1}}_q$$
$$ \times \kbrack{\mkern-12mu\mbox{
    \tikz[baseline=-.6ex, scale=.08]{
        \setsquarecoords{12}
         
        \begin{scope}
             
            \draw[warc] (-8,-6) to[out=east, in=west] (0,6);
            \draw[warc] (-8,6) to[out=east, in=west] (0,-6);
            \draw[warc] (0,-6) to[out=east, in=west] (8,6);
            \draw[warc] (0,6) to[out=east, in=west] (8,-6);

            \draw[warc] (-10, 4) -- (-8, 4);
            \draw[warc] (-8, 4) to[out=east, in=east] (-8, -4);
            \draw[warc] (-10, -4) -- (-8, -4);
            \node at (-13,0) [align=center] {\tiny $m-k_i$};
             
            \draw[warc] (-14, 6) -- (-8, 6);
            
            \draw[warc] (-14, -6) -- (-8, -6);

            \node[vJW, scale=0.06, label={[xshift=2pt, yshift=6pt]left:\scriptsize $n$}] at (-10,6) {};
            
            \node[vJW, scale=0.06, label={[xshift=2pt, yshift=-6pt]left:\scriptsize $m$}] at (-10,-6) {};
            
            \node at (20,0) [align=center] {\scriptsize $...$};
            \node at (20,0) [above, align=center] {\scriptsize $p-i-1$\\[-8pt] \scriptsize full twists};

            \node at (0.5,-5) [below, align=center] {\tiny $k_i + n-m$};
            \node at (0.5,5) [above, align=center] {\tiny $k_i$};
        \end{scope}

        \begin{scope}[shift={(40,0)}]
            \draw[warc] (-8,-6) to[out=east, in=west] (0,6);
            \draw[warc] (-8,6) to[out=east, in=west] (0,-6);
            \draw[warc] (0,-6) to[out=east, in=west] (8,6);
            \draw[warc] (0,6) to[out=east, in=west] (8,-6);

            \draw[warc] (8, 6) -- (14, 6);
            \draw[warc] (8, -6) -- (14, -6);

            \draw[warc] (10, 4) -- (18, 4);
            \draw[warc] (18, 4) to[out=east, in=east] (18, -4);
            \draw[warc] (10, -4) -- (18, -4);
            \node at (12.5,0) [align=center] {\tiny $k_i-k_{i+1}$};

            \draw[warc] (32, 4) -- (34, 4);
            \draw[warc] (32, 4) to[out=west, in=west] (32, -4);
            \draw[warc] (32, -4) -- (34, -4);
            \node at (30,0) [align=center, right] {\tiny $m-k_{i+1}$};

            \node at (22,5) [above, align=center] {\tiny $k_{i+1} + n - m$};
            \node at (24,-5) [below, align=center] {\tiny $k_{i+1}$};

            \node[vJW, scale=0.06] at (10,-6) {};
            \node[vJW, scale=0.06] at (10,6) {};

            \node at (-0.5,-5) [below, align=center] {\tiny $k_i + n-m$};
            \node at (0.5,5) [above, align=center] {\tiny $k_i$};

            \draw[warc] (14, 6) -- (38, 6);
            \draw[warc] (14, -6) -- (38, -6);
            \node[vJW, scale=0.06, label={[xshift=-2pt, yshift=-6pt]right:\scriptsize $m$}] at (34,-6) {};
            \node[vJW, scale=0.06, label={[xshift=-2pt, yshift=6pt]right:\scriptsize $n$}] at (34,6) {};
        \end{scope}
    }
}\mkern-12mu}$$

$$ = \sum_{k_{i+1}=0}^{k_i} q^{\frac12(2k_{i+1}^2 + 2k_{i+1}(n-m) - k_i(k_i + n - m) - k_i + k_{i+1})} (-1)^{k_i - k_{i+1}} \frac{(q)_{k_i + n - m}}{(q)_{k_{i+1} + n - m}}\binom{k_i}{k_i - k_{i+1}}_q$$

$$\times \left(\frac{\mu_{k_{i+1}}\mu_{k_{i+1}+n - m}}{\mu_{k_i}\mu_{k_i + n - m}}\right)^{p-i-1} \kbrack{\sigma^{(p-i-1)}_{m,n}(k_{i+1})}$$

Then, we have a recurrence relation for $\kbrack{\sigma_{m,n}^{(p-i)}(k_{i})}$ in terms of $\kbrack{\sigma_{m,n}^{(p-i-1)}(k_{i+1})}$, which allows us to express $\kbrack{\sigma_{m,n}^{(p)}(k_{0})}$ in terms of $\kbrack{\sigma_{m,n}^{(0)}(k_{p})} $where $k_0 = m$, with coefficients depending on $k_i$ for $i = 0, 1, \cdots, p-1$. Let us use the label $k'_{i+1} = k_i - k_{i+1}$. Many of the terms end up telescoping and simplifying:

$$
\kbrack{\sigma^{(p)}_{m,n}(k_0)} = \sum_{\substack{0 \leq k_p \leq \cdots \\ \leq k_1 \leq k_0 = m}} (-1)^{m - k_p} \cdot \left(\prod_{i=0}^{p-1}q^{\frac12(2k_{i+1}^2 - k_i^2)}q^{\frac{(n-m)}{2}(2k_{i+1}-k_i)}q^{\frac12(k_{i+1}-k_i)}\right) \cdot \frac{(q)_{n}}{(q)_{k_p + n - m}}
$$
$$
\times \binom{m}{k_1', k_2', \cdots, k_{p}', k_p}_q \times q^{\frac12(p-1)(-2m-mn)}q^{\frac12\sum_{i=1}^{p-1}[k_i^2 + k_i(n-m+2)]} \times \kbrack{\sigma^{0}_{m,n}(k_p)}$$
We can simplify the algebra by considering the total exponent $E$ from the $q$ term and $F$ from the $\mu$ term, then verifying this matches the $m=n$ case of Yuasa.
$$E= \frac12\left(-k_0^2 + 2k_p^2 + \sum_{i=1}^{p-1}k_i^2\right) + \frac{(n-m)}{2}\left(-k_0 + 2k_p + \sum_{i=1}^{p-1}k_i\right) + \frac{-k_0 + k_p}{2}$$
$$F = \frac{1}{2}\left((p-1)(-2m-mn)  + \sum_{i=1}^{p-1}[k_i^2 + k_i(n-m+2)]\right)$$
$$E+F = \frac12\left(-m^2 + 2k_p^2 + (n-m)(2k_p - m) - m + k_p +(p-1)(-2m-mn)\right)$$
$$ + \frac12 \sum_{i=1}^{p-1}[2k_i^2 + 2k_i(n-m+1)]$$
$$ \implies E+F = \frac12\left(2k_p(n-m) + 2k_p^2 + m + k_p + p(-2m-mn)\right) + \sum_{i=1}^{p-1}k_i^2 + k_i(n-m+1)$$
$$ = \frac12\left(m - k_p + p(-2m-mn)\right) +\sum^p_{i=1}[k_i^2 + k_i(n-m+1)]$$
This gives the exponent we were seeking in the above formula. One can see that setting $m=n$ gives us the formula found in Proposition 5.4. \hfill$\square$

\section{Application}
\textbf{Definition 7.1 $\llbracket$2-Bridge Links$\rrbracket$:} One useful place to apply our theory for full twists of differently colored strands is in the 2-colored Jones polynomials of a family of links known as \textit{2-bridge links}. These are composed of series of full twists, labeled $L =[2a_1, 2a_2, \cdots 2a_l]$ with $a_1, \cdots, a_l \in \zz_{\neq 0}$. Negative $a_i$ correspond to negatively oriented twists and are associated with a swapping in the full-twist formula of $q \ra q^{-1}$.
$$\mbox{
    \text{for $l$ odd:} \quad \tikz[baseline=-.6ex, scale=.06]{
        \setsquarecoords{12}
         
        \begin{scope}
             
            \draw[warc] (-50, 12) -- (50, 12);
            \draw[warc] (-50, 4) -- (50, 4);
            \draw[warc] (-50, -12) -- (50, -12);
            \draw[warc] (-50, -4) -- (50, -4);

            \draw[warc] (-50, 12) to[out=west, in=west] (-50, 4);
            \draw[warc] (-50, -12) to[out=west, in=west] (-50, -4);
            \draw[warc] (50, 12) to[out=east, in=east] (50, 4);
            \draw[warc] (50, -12) to[out=east, in=east] (50, -4);

            \draw[fill=white] (-44,-6) rectangle (-32,6) node[midway] {$2a_1$};
            \draw[fill=white] (-28,-14) rectangle (-16,-2) node[midway] {$2a_2$};
            \draw[fill=white] (-12,-6) rectangle (0,6) node[midway] {$2a_3$};
            \draw[fill=white] (4,-14) rectangle (16,-2) node[midway] {$2a_4$};
            \draw[fill=white] (32,-6) rectangle (44,6) node[midway] {$2a_l$};

            \draw[fill=white, draw=none] (18,-13) rectangle (28,7) node[midway] {$\substack{\cdots\\ \cdots}$};

        \end{scope}
    }
}$$
$$\mbox{
    \text{for $l$ even:} \quad \tikz[baseline=-.6ex, scale=.06]{
        \setsquarecoords{12}
         
        \begin{scope}
             
            \draw[warc] (-50, 12) -- (50, 12);
            \draw[warc] (-50, 4) -- (50, 4);
            \draw[warc] (-50, -12) -- (50, -12);
            \draw[warc] (-50, -4) -- (50, -4);

            \draw[warc] (-50, 12) to[out=west, in=west] (-50, 4);
            \draw[warc] (-50, -12) to[out=west, in=west] (-50, -4);
            \draw[warc] (50, -4) to[out=east, in=east] (50, 4);
            \draw[warc] (50, -12) to[out=east, in=east] (50, 12);

            \draw[fill=white] (-44,-6) rectangle (-32,6) node[midway] {$2a_1$};
            \draw[fill=white] (-28,-14) rectangle (-16,-2) node[midway] {$2a_2$};
            \draw[fill=white] (-12,-6) rectangle (0,6) node[midway] {$2a_3$};
            \draw[fill=white] (4,-14) rectangle (16,-2) node[midway] {$2a_4$};
            \draw[fill=white] (32,-14) rectangle (44,-2) node[midway] {$2a_l$};

            \draw[fill=white, draw=none] (18,-13) rectangle (28,7) node[midway] {$\substack{\cdots\\ \cdots}$};

        \end{scope}
    }
}$$

Yuasa explored $\spl_2$ and $\spl_3$ colored Jones polynomials for 2-bridge links of the same color in \cite{yuasa}. In this paper we generalize the formula for differently colored links. Note, that in the case of $l$ being even, the diagram will actually result in having only one single-colored component. For the purposes of using our generalized formula, we will focus on the case where $l$ is odd, and the two components separate and allow us to color them with $m$ and $n$. Let us begin to generalize using the same proof method of Yuasa to find the 2-colored Jones polynomial for the odd 2-bridge links. \

\hspace{-0.25in}\textbf{Definition 7.2 $\llbracket$2-Color Jones Invariant$\rrbracket$:}
The colored Jones polynomials are usually defined in terms of one color, but there are analogues which take into account multiple strands of different colors. In our case, two will suffice. Let $L_{m,n} = L_m \cup L_n$, the two components of a link which are colored by $m,n \in \zz^+$ respectively, and $w(L_i)$ is the writhe of the component $L_i$. Then the \textit{2-color Jones polynomials} are given by:
$$J_{m,n}(L_{m,n}, q) = \mu_m ^{-w(L_m)} \mu_n^{-w(L_n)} \kbrack{L_{m,n}}$$

\hspace{-0.25in}\textbf{Lemma 7.3}
$$\Biggl \langle \mkern-8mu
\mbox{
    \tikz[baseline=-.6ex, scale=.08]{
        \setsquarecoords{12}
         
        \begin{scope}
             
            \draw[warc] (-8,8) -- (8,8) node[midway, below, black, yshift=2pt] {\scriptsize $n-k$};
            \draw[warc] (-8,-8) -- (8,-8) node[midway, above, black, yshift=-2pt] {\scriptsize $m-k$};

            \draw[warc] (-8, 4) to[out=east, in =east, looseness=1.5, edge node = {node[ left, midway, black, font=\scriptsize] {$k$}}] (-8, -4);
             
            \draw[warc] (8, 4) to[out=west, in =west, looseness=1.5, edge node = {node[ right, midway, black, font=\scriptsize] {$k$}}] (8, -4);

            \draw[warc] (-12, 6) -- (-8, 6) node[midway, above, black, yshift=0pt, xshift=-1pt]{\scriptsize $n$};
            
            \draw[warc] (-12, -6) -- (-8, -6) node[midway, below, black, yshift=0pt, xshift=-2pt]{\scriptsize $m$};
            
            \draw[warc] (12, 6) -- (8, 6) node[midway, above, black, yshift=0pt, xshift=1pt]{\scriptsize $n$};
            
            \draw[warc] (12, -6) -- (8, -6) node[midway, below, black, yshift=0pt, xshift=2pt]{\scriptsize $m$};

            \draw[warc] (12, -6) to[out=east, in=east] (12, -12);
            \draw[warc] (-12, -6) to[out=west, in=west] (-12, -12);
            \draw[warc] (-12, -12) -- (12, -12);

            \draw[warc] (12, 6) to[out=east, in=east] (12, 12);
            \draw[warc] (-12, 6) to[out=west, in=west] (-12, 12);
            \draw[warc] (-12, 12) -- (12, 12);
             
            \node[vJW, scale=0.06] at (-8,6) {};
            \node[vJW, scale=0.06] at (-8,-6) {};
            \node[vJW, scale=0.06] at (8,6) {};
            \node[vJW, scale=0.06] at (8,-6) {};
        \end{scope}
    }
} \mkern-8mu \Biggr \rangle = (-1)^{n+m-k}\frac{[n+1][m+1]}{[k+1]}$$
This follows directly from the properties in Propositions 4.3 and 4.5, and helps in the computation of the colored Jones polynomials of the odd 2-bridge links.\

\hspace{-0.25in}\textbf{Notation:} For the sake of cleaning up the formulas used in the next theorem. Let us define $K_0 = n$ and $K_j = m  - k_{|a_j|}^{(j)}$, with $k^{(j)}_0 = K_j$ and $k_{i+1}^{(j+1)'} = k_i^{(j)}-k_{i+1}^{(j)}$.

$$\text{Twist Sign: } \epsilon_{j+1} = \frac{a_{j+1}}{|a_{j+1}|}$$
$$\text{Sliding Formula: } \delta_{n,m}(K_j) = \left(\frac{\mu_{K_j}\mu_{K_j+n-m}}{\mu_n\mu_m}\right)^{\epsilon_{j+1}}$$
$$\text{Full Twist Coefficient: }\gamma_{m,n}(k_1^{(j+1)}, \cdots, k^{(j+1)}_{|a_{j+1}|}; q) = (-1)^{m - k_{|a_{j+1}|}^{(j+1)}}q^{-\frac{|a_{j+1}|}{2}(mn + 2m)}q^{\frac{1}{2}(m - k_{|a_{j+1}|}^{(j+1)})}$$
$$\times q^{\sum_{i=1}^{|a_{j+1}|}(k_i^{(j+1)})^2 + (k_{i}^{(j+1)})(n-m+1)} \frac{(q)_{n}}{(q)_{k_{|a_{j+1}|}^{(j+1)}+n - m}} \binom{m}{k_1^{(j+1)'}, k_2^{(j+1)'}, \cdots, k_{|a_{j+1}|}^{(j+1)'}, k_{|a_{j+1}|}^{(j+1)}}$$

\hspace{-0.25in}\textbf{Theorem 7.4 $\llbracket$2-Colored Jones Polynomial of Odd $(m,n)$-Colored 2-Bridge Links$\rrbracket$:} For $m,n \in \zz^+$, $m \leq n$ and $L=[2a_1, \cdots, 2a_l]$ with $l$ odd, let us define ${L_{m,n}}$ as the coloring of $L$ with the top strand replaced by $n$ and the bottom strand replaced by $m$ as such:
$$L_{m,n} \coloneq \quad \mbox{
    \tikz[baseline=-.6ex, scale=.06]{
        \setsquarecoords{12}
         
        \begin{scope}
             
            \draw[warc] (-55, 12) -- (50, 12);
            \draw[warc] (-55, 4) -- (50, 4);
            \draw[warc] (-55, -12) -- (50, -12);
            \draw[warc] (-55, -4) -- (50, -4);

            \draw[warc] (-55, 12) to[out=west, in=west] (-55, 4);
            \draw[warc] (-55, -12) to[out=west, in=west] (-55, -4);
            \draw[warc] (50, 12) to[out=east, in=east] (50, 4);
            \draw[warc] (50, -12) to[out=east, in=east] (50, -4);

            \draw[fill=white] (-44,-6) rectangle (-32,6) node[midway] {$2a_1$};
            \draw[fill=white] (-28,-14) rectangle (-16,-2) node[midway] {$2a_2$};
            \draw[fill=white] (-12,-6) rectangle (0,6) node[midway] {$2a_3$};
            \draw[fill=white] (4,-14) rectangle (16,-2) node[midway] {$2a_4$};
            \draw[fill=white] (32,-6) rectangle (44,6) node[midway] {$2a_l$};

            \draw[fill=white, draw=none] (18,-13) rectangle (28,7) node[midway] {$\substack{\cdots\\ \cdots}$};
            
            \node[vJW, scale=0.08] at (-50,-12) {};
            \node[vJW, scale=0.08] at (-50,12) {};

            \node at (-46,-20) [above, align=center] {\scriptsize $m$};
            \node at (-46,20) [below, align=center] {\scriptsize $n$};
        \end{scope}
    }
}\quad \text{with $l$ odd}$$
Then: $$J_{m,n}(L_{m,n}, q) =$$$$\left(\prod_{\substack{0 \leq j \leq l -1 \\ \text{$j$ even}}}\sum_{0 \leq k^{(j+1)}_{|a_{j+1}|} \leq \cdots \leq k^{(j+1)}_{1}\leq K_j}\delta_{K_j, K_j+n-m}(m-K_j)\ \gamma_{{K_j, K_j + n - m}}\left(k^{(j+1)}_{1}, \cdots , k^{(j+1)}_{|a_{j+1}|}; q^{\epsilon_{j+1}}\right)\right)$$
$$\times \left(\prod_{\substack{0 \leq j \leq l -1 \\ \text{$j$ odd}}}\sum_{0 \leq k^{(j+1)}_{|a_{j+1}|} \leq \cdots \leq k^{(j+1)}_{1}\leq K_j}\delta_{m,m}(m-K_j)\ \gamma_{{K_j, K_j}}\left(k^{(j+1)}_{1}, \cdots , k^{(j+1)}_{|a_{j+1}|}; q^{\epsilon_{j+1}}\right)\right)$$
$$\times \mu_m^{-2(a_2 + a_4 + \cdots + a_{l-1})}\mu_n^{-2(a_1 + a_3 + \cdots + a_{l})}(-1)^{n+K_l}\frac{[n+1][m+1]}{[m-K_l + 1]}$$

\hspace{-0.25in}\textbf{Proof:} We have already defined the full twist formula for $(m,n)$-colored strands as well as the sliding formula, so we are able to reduce everything individually using the colored Kauffman bracket. Let us reduce from left to right, calculating the Kauffman bracket of the following diagrams (from which the 2-bridge link is composed)\

\hspace{-0.25in}\textit{Reduction for middle twists:}
$$\Biggl \langle \mkern-8mu \mbox{\tikz[baseline=-.6ex, scale=.08]{
        \setsquarecoords{12}
         
        \begin{scope}
             
            \draw[warc] (-20, 12) -- (20, 12);
            \draw[warc] (-20, 4) -- (20, 4);
            \draw[warc] (-20, -12) -- (20, -12);
            \draw[warc] (-20, -4) -- (20, -4);

            \draw[warc] (-20, 11) to[out=west, in=west, looseness=1.5] (-20, 5);
            \draw[warc] (-20, 3) to[out=west, in=west, looseness=1.5] (-20, -3);
            \draw[warc] (-20, -5) to[out=west, in=west, looseness=1.5] (-20, -11);

            \draw[warc] (-20, 13) -- (-44, 13);
            \draw[warc] (-44, 13) to[out=west, in=west, looseness=0.75] (-44, -13);
            \draw[warc] (-20, -13) -- (-44, -13);

            \node[vJW, scale=0.04] at (-20,12) {};
            \node[vJW, scale=0.04] at (-20,4) {};
            \node[vJW, scale=0.04] at (-20,-12) {};
            \node[vJW, scale=0.04] at (-20,-4) {};
            \node[vJW, scale=0.04] at (20,12) {};
            \node[vJW, scale=0.04] at (20,4) {};
            \node[vJW, scale=0.04] at (20,-4) {};
            \node[vJW, scale=0.04] at (20,-12) {};
            
            \node at (-21,0) [left, align=center] {\scriptsize $m - K_j$};
            \node at (-21,8) [left, align=center] {\scriptsize $n-m + K_j$};
            \node at (-21,-8) [left, align=center] {\scriptsize $K_j$};

            \node at (-50,0) [left, align=center] {\scriptsize $m - K_j$};
            
            \node at (14,11) [above, align=center] {\scriptsize $n$};
            \node at (14,3) [above, align=center] {\scriptsize $n$};
            \node at (14,-5) [above, align=center] {\scriptsize $m$};
            \node at (14,-13) [above, align=center] {\scriptsize $m$};

            \node at (-14,11) [above, align=center] {\scriptsize $n$};
            \node at (-14,3) [above, align=center] {\scriptsize $n$};
            \node at (-14,-5) [above, align=center] {\scriptsize $m$};
            \node at (-14,-13) [above, align=center] {\scriptsize $m$};

            \draw[fill=white] (-8,-8) rectangle (8,8) node[midway] {\scriptsize $2a_{j+1}$};
        \end{scope}
    }
} \Biggr \rangle $$
$$= \delta_{n,m}(m-K_j)^{2|a_{j+1}|} \Biggl \langle \mkern-8mu \mbox{\tikz[baseline=-.6ex, scale=.08]{
        \setsquarecoords{12}
         
        \begin{scope}
             
            \draw[warc] (-20, 12) -- (20, 12);
            \draw[warc] (-20, 4) -- (20, 4);
            \draw[warc] (-20, -12) -- (20, -12);
            \draw[warc] (-20, -4) -- (20, -4);

            \draw[warc] (-20, 11) to[out=west, in=west, looseness=1.5] (-20, 4);
            \draw[warc] (20, 3) to[out=west, in=west, looseness=1.5] (20, -3);
            \draw[warc] (-20, -4) to[out=west, in=west, looseness=1.5] (-20, -11);

            \draw[warc] (-20, 13) -- (-44, 13);
            \draw[warc] (-44, 13) to[out=west, in=west, looseness=0.75] (-44, -13);
            \draw[warc] (-20, -13) -- (-44, -13);

            \node[vJW, scale=0.04] at (-20,12) {};
            \node[vJW, scale=0.04] at (-20,-12) {};
            \node[vJW, scale=0.04] at (20,12) {};
            \node[vJW, scale=0.04] at (20,4) {};
            \node[vJW, scale=0.04] at (20,-4) {};
            \node[vJW, scale=0.04] at (20,-12) {};
            
            \node at (19,0) [right, align=center] {\scriptsize $m - K_j$};
            \node at (-21,8) [left, align=center] {\scriptsize $n-m + K_j$};
            \node at (-21,-8) [left, align=center] {\scriptsize $K_j$};

            \node at (-50,0) [left, align=center] {\scriptsize $m - K_j$};
            
            \node at (14,11) [above, align=center] {\scriptsize $n$};
            \node at (14,-13) [above, align=center] {\scriptsize $m$};

            \node at (-14,11) [above, align=center] {\scriptsize $n$};
            \node at (-14,-13) [above, align=center] {\scriptsize $m$};

            \draw[fill=white] (-8,-8) rectangle (8,8) node[midway] {\scriptsize $2a_{j+1}$};
        \end{scope}
    }
} \Biggr \rangle $$
$$= \delta_{n,m}(m-K_j)^{2|a_{j+1}|} \sum_{\substack{0 \leq k^{(j+1)}_{|a_{j+1}|} \leq \cdots \leq k^{(j+1)}_{1} \leq K_j}} \gamma_{K_j, K_j + n - m}(k_1^{(j+1)}, \cdots, k_{|a_{j+1}|}^{(j+1)}; q^{\epsilon_{j+1}})$$
$$\times \Biggl \langle \mkern-8mu \mbox{\tikz[baseline=-.6ex, scale=.08]{
        \setsquarecoords{12}
         
        \begin{scope}
             
            \draw[warc] (-20, 12) -- (20, 12);
            \draw[warc] (-20, 4) -- (20, 4);
            \draw[warc] (-20, -12) -- (20, -12);
            \draw[warc] (-20, -4) -- (20, -4);

            \draw[warc] (-20, 11) to[out=west, in=west, looseness=1.5] (-20, 4);
            \draw[warc] (20, 3) to[out=west, in=west, looseness=1.5] (20, -3);
            \draw[warc] (-20, -4) to[out=west, in=west, looseness=1.5] (-20, -11);

            \draw[warc] (-20, 13) -- (-56, 13);
            \draw[warc] (-56, 13) to[out=west, in=west, looseness=0.75] (-56, -13);
            \draw[warc] (-20, -13) -- (-56, -13);

            \node[vJW, scale=0.04] at (-20,12) {};
            \node[vJW, scale=0.04] at (-20,-12) {};
            \node[vJW, scale=0.04] at (20,12) {};
            \node[vJW, scale=0.04] at (20,4) {};
            \node[vJW, scale=0.04] at (20,-4) {};
            \node[vJW, scale=0.04] at (20,-12) {};
            
            \node at (19,0) [right, align=center] {\scriptsize $m - k_{|a_{j+1}|}^{(j+1)}$};
            \node at (-21,8) [left, align=center] {\scriptsize $n-m + k_{|a_{j+1}|}^{(j+1)}$};
            \node at (-21,-8) [left, align=center] {\scriptsize $k_{|a_{j+1}|}^{(j+1)}$};

            \node at (-60,0) [left, align=center] {\scriptsize $m - k_{|a_{j+1}|}^{(j+1)}$};
            
            \node at (14,11) [above, align=center] {\scriptsize $n$};
            \node at (14,-13) [above, align=center] {\scriptsize $m$};

            \node at (-14,11) [above, align=center] {\scriptsize $n$};
            \node at (-14,-13) [above, align=center] {\scriptsize $m$};
        \end{scope}
    }
} \Biggr \rangle $$
$$= \delta_{n,m}(m-K_j)^{2|a_{j+1}|} \sum_{\substack{0 \leq k^{(j+1)}_{|a_{j+1}|} \leq \cdots \leq k^{(j+1)}_{1} \leq K_j}} \gamma_{K_j, K_j + n - m}(k_1^{(j+1)}, \cdots, k_{|a_{j+1}|}^{(j+1)}; q^{\epsilon_{j+1}})$$
$$\times \Biggl \langle \mkern-8mu \mbox{\tikz[baseline=-.6ex, scale=.08]{
        \setsquarecoords{12}
         
        \begin{scope}
             
            \draw[warc] (20, 11) to[out=west, in=west, looseness=1.5] (20, 4);
            \draw[warc] (20, 3) to[out=west, in=west, looseness=1.5] (20, -3);
            \draw[warc] (20, -4) to[out=west, in=west, looseness=1.5] (20, -11);

            \draw[warc] (-5, 12) -- (20, 12);
            \draw[warc] (-5, 12) to[out=west, in=west, looseness=0.75] (-5, -12);
            \draw[warc] (-5, -12) -- (20, -12);

            \node[vJW, scale=0.04] at (20,12) {};
            \node[vJW, scale=0.04] at (20,4) {};
            \node[vJW, scale=0.04] at (20,-4) {};
            \node[vJW, scale=0.04] at (20,-12) {};
            
            \node at (19,0) [left, align=center] {\scriptsize $K_{j+1}$};
            \node at (19,8) [left, align=center] {\scriptsize $n- K_{j+1}$};
            \node at (19,-8) [left, align=center] {\scriptsize $m-K_j$};

            \node at (-10,0) [left, align=center] {\scriptsize $K_{j+1}$};
            
        \end{scope}
    }
} \Biggr \rangle $$

\hspace{-0.25in}\textit{Reduction for lower twists:}
$$\Biggl \langle \mkern-8mu \mbox{\tikz[baseline=-.6ex, scale=.08]{
        \setsquarecoords{12}
         
        \begin{scope}
             
            \draw[warc] (-20, 12) -- (20, 12);
            \draw[warc] (-20, 4) -- (20, 4);
            \draw[warc] (-20, -12) -- (20, -12);
            \draw[warc] (-20, -4) -- (20, -4);

            \draw[warc] (-20, 11) to[out=west, in=west, looseness=1.5] (-20, 5);
            \draw[warc] (-20, 3) to[out=west, in=west, looseness=1.5] (-20, -3);
            \draw[warc] (-20, -5) to[out=west, in=west, looseness=1.5] (-20, -11);

            \draw[warc] (-20, 13) -- (-44, 13);
            \draw[warc] (-44, 13) to[out=west, in=west, looseness=0.75] (-44, -13);
            \draw[warc] (-20, -13) -- (-44, -13);

            \node[vJW, scale=0.04] at (-20,12) {};
            \node[vJW, scale=0.04] at (-20,4) {};
            \node[vJW, scale=0.04] at (-20,-12) {};
            \node[vJW, scale=0.04] at (-20,-4) {};
            \node[vJW, scale=0.04] at (20,12) {};
            \node[vJW, scale=0.04] at (20,4) {};
            \node[vJW, scale=0.04] at (20,-4) {};
            \node[vJW, scale=0.04] at (20,-12) {};
            
            \node at (-21,0) [left, align=center] {\scriptsize $K_j$};
            \node at (-21,8) [left, align=center] {\scriptsize $n-K_{j}$};
            \node at (-21,-8) [left, align=center] {\scriptsize $m-K_j$};

            \node at (-50,0) [left, align=center] {\scriptsize $K_j$};
            
            \node at (14,11) [above, align=center] {\scriptsize $n$};
            \node at (14,3) [above, align=center] {\scriptsize $n$};
            \node at (14,-5) [above, align=center] {\scriptsize $m$};
            \node at (14,-13) [above, align=center] {\scriptsize $m$};

            \node at (-14,11) [above, align=center] {\scriptsize $n$};
            \node at (-14,3) [above, align=center] {\scriptsize $n$};
            \node at (-14,-5) [above, align=center] {\scriptsize $m$};
            \node at (-14,-13) [above, align=center] {\scriptsize $m$};

            \draw[fill=white] (-8,-16) rectangle (8,0) node[midway] {\scriptsize $2a_{j+1}$};
        \end{scope}
    }
} \Biggr \rangle $$
$$= \delta_{m,m}(m-K_j)^{2|a_{j+1}|} \Biggl \langle \mkern-8mu \mbox{\tikz[baseline=-.6ex, scale=.08]{
        \setsquarecoords{12}
         
        \begin{scope}
             
            \draw[warc] (-20, 12) -- (20, 12);
            \draw[warc] (-20, 4) -- (20, 4);
            \draw[warc] (-20, -12) -- (20, -12);
            \draw[warc] (-20, -4) -- (20, -4);

            \draw[warc] (-20, 11) to[out=west, in=west, looseness=1.5] (-20, 5);
            \draw[warc] (-20, 3) to[out=west, in=west, looseness=1.5] (-20, -4);
            \draw[warc] (20, -5) to[out=west, in=west, looseness=1.5] (20, -11);

            \draw[warc] (-20, 13) -- (-44, 13);
            \draw[warc] (-44, 13) to[out=west, in=west, looseness=0.75] (-44, -12);
            \draw[warc] (-20, -12) -- (-44, -12);

            \node[vJW, scale=0.04] at (-20,12) {};
            \node[vJW, scale=0.04] at (-20,4) {};
            \node[vJW, scale=0.04] at (20,12) {};
            \node[vJW, scale=0.04] at (20,4) {};
            \node[vJW, scale=0.04] at (20,-4) {};
            \node[vJW, scale=0.04] at (20,-12) {};
            
            \node at (-21,0) [left, align=center] {\scriptsize $K_j$};
            \node at (-21,8) [left, align=center] {\scriptsize $n-K_{j}$};
            \node at (21,-8) [right, align=center] {\scriptsize $m-K_j$};

            \node at (-50,0) [left, align=center] {\scriptsize $K_j$};
            
            \node at (14,11) [above, align=center] {\scriptsize $n$};
            \node at (14,3) [above, align=center] {\scriptsize $n$};

            \node at (-14,11) [above, align=center] {\scriptsize $n$};
            \node at (-14,3) [above, align=center] {\scriptsize $n$};

            \draw[fill=white] (-8,-16) rectangle (8,0) node[midway] {\scriptsize $2a_{j+1}$};
        \end{scope}
    }
} \Biggr \rangle $$
$$= \delta_{m,m}(m-K_j)^{2|a_{j+1}|} \sum_{\substack{0 \leq k^{(j+1)}_{|a_{j+1}|} \leq \cdots \leq k^{(j+1)}_{1} \leq K_j}} \gamma_{K_j, K_j}(k_1^{(j+1)}, \cdots, k_{|a_{j+1}|}^{(j+1)}; q^{\epsilon_{j+1}})$$
$$\times \Biggl \langle \mkern-8mu \mbox{\tikz[baseline=-.6ex, scale=.08]{
        \setsquarecoords{12}
         
        \begin{scope}
             
            \draw[warc] (-20, 12) -- (20, 12);
            \draw[warc] (-20, 4) -- (20, 4);
            \draw[warc] (-20, -12) -- (20, -12);
            \draw[warc] (-20, -4) -- (20, -4);

            \draw[warc] (-20, 11) to[out=west, in=west, looseness=1.5] (-20, 5);
            \draw[warc] (-20, 3) to[out=west, in=west, looseness=1.5] (-20, -4);
            \draw[warc] (20, -5) to[out=west, in=west, looseness=1.5] (20, -11);

            \draw[warc] (-20, -5) to[out=east, in=east, looseness=1.5] (-20, -11);

            \draw[warc] (-20, 13) -- (-44, 13);
            \draw[warc] (-44, 13) to[out=west, in=west, looseness=0.75] (-44, -12);
            \draw[warc] (-20, -12) -- (-44, -12);

            \node[vJW, scale=0.04] at (-20,12) {};
            \node[vJW, scale=0.04] at (-20,4) {};
            \node[vJW, scale=0.04] at (20,12) {};
            \node[vJW, scale=0.04] at (20,4) {};
            \node[vJW, scale=0.04] at (20,-4) {};
            \node[vJW, scale=0.04] at (20,-12) {};
            \node[vJW, scale=0.04] at (-20,-4) {};
            \node[vJW, scale=0.04] at (-20,-12) {};
            
            \node at (-21,0) [left, align=center] {\scriptsize $K_j$};
            \node at (-21,8) [left, align=center] {\scriptsize $n-K_{j}$};
            \node at (21,-8) [right, align=center] {\scriptsize $m-k_{|a_{j+1}|}^{(j+1)}$};
            \node at (-18,-8) [right, align=center] {\scriptsize $K_j - k_{|a_{j+1}|}^{(j+1)}$};
             \node at (-10,0) [right, align=center] {\scriptsize $k_{|a_{j+1}|}^{(j+1)}$};
             \node at (-10,-18) [right, align=center] {\scriptsize $k_{|a_{j+1}|}^{(j+1)}$};

            \node at (-50,0) [left, align=center] {\scriptsize $K_j$};
            
            \node at (14,11) [above, align=center] {\scriptsize $n$};
            \node at (14,3) [above, align=center] {\scriptsize $n$};

            \node at (-14,11) [above, align=center] {\scriptsize $n$};
            \node at (-14,3) [above, align=center] {\scriptsize $n$};
        \end{scope}
    }
} \Biggr \rangle $$
$$= \delta_{m,m}(m-K_j)^{2|a_{j+1}|} \sum_{\substack{0 \leq k^{(j+1)}_{|a_{j+1}|} \leq \cdots \leq k^{(j+1)}_{1} \leq K_j}} \gamma_{K_j, K_j}(k_1^{(j+1)}, \cdots, k_{|a_{j+1}|}^{(j+1)}; q^{\epsilon_{j+1}})$$
$$\times \Biggl \langle \mkern-8mu \mbox{\tikz[baseline=-.6ex, scale=.08]{
        \setsquarecoords{12}
         
        \begin{scope}
             
            \draw[warc] (20, 11) to[out=west, in=west, looseness=1.5] (20, 4);
            \draw[warc] (20, 3) to[out=west, in=west, looseness=1.5] (20, -3);
            \draw[warc] (20, -4) to[out=west, in=west, looseness=1.5] (20, -11);

            \draw[warc] (-10, 12) -- (20, 12);
            \draw[warc] (-10, 12) to[out=west, in=west, looseness=0.75] (-10, -12);
            \draw[warc] (-10, -12) -- (20, -12);

            \node[vJW, scale=0.04] at (20,12) {};
            \node[vJW, scale=0.04] at (20,4) {};
            \node[vJW, scale=0.04] at (20,-4) {};
            \node[vJW, scale=0.04] at (20,-12) {};
            
            \node at (19,0) [left, align=center] {\scriptsize $m -K_{j+1}$};
            \node at (19,8) [left, align=center] {\scriptsize $n- m +K_{j+1}$};
            \node at (19,-8) [left, align=center] {\scriptsize $K_{j+1}$};

            \node at (-14,0) [left, align=center] {\scriptsize $m - K_{j+1}$};
            
        \end{scope}
    }
} \Biggr \rangle $$
Then in computing the bracket of the entire link, we should handle $a_j$'s separately for $j$ odd versus $j$ even:
$$\kbrack{L_{m,n}} =$$$$ \left(\prod_{\substack{0 \leq j \leq l -1 \\ \text{$j$ even}}}\sum_{0 \leq k^{(j+1)}_{|a_{j+1}|} \leq \cdots \leq k^{(j+1)}_{1}\leq K_j}\delta_{n,m}(m-K_j)\ \gamma_{{K_j, K_j + n - m}}\left(k^{(j+1)}_{1}, \cdots , k^{(j+1)}_{|a_{j+1}|}; q^{\epsilon_{j+1}}\right)\right)$$
$$\times \left(\prod_{\substack{0 \leq j \leq l -1 \\ \text{$j$ odd}}}\sum_{0 \leq k^{(j+1)}_{|a_{j+1}|} \leq \cdots \leq k^{(j+1)}_{1}\leq K_j}\delta_{m,m}(m-K_j)\ \gamma_{{K_j, K_j}}\left(k^{(j+1)}_{1}, \cdots , k^{(j+1)}_{|a_{j+1}|}; q^{\epsilon_{j+1}}\right)\right)$$
$$\times \Biggl \langle \mkern-8mu
\mbox{
    \tikz[baseline=-.6ex, scale=.08]{
        \setsquarecoords{12}
         
        \begin{scope}
             
            \draw[warc] (-8,8) -- (14,8) node[midway, below, black, yshift=2pt] {\scriptsize $K_l + n -m$};
            \draw[warc] (-8,-8) -- (14,-8) node[midway, above, black, yshift=-2pt] {\scriptsize $m-K_l$};

            \draw[warc] (-8, 4) to[out=east, in =east, looseness=1.5, edge node = {node[ left, midway, black, font=\scriptsize] {$m-K_l$}}] (-8, -4);
             
            \draw[warc] (14, 4) to[out=west, in =west, looseness=1.5, edge node = {node[ right, midway, black, font=\scriptsize] {$m-K_l$}}] (14, -4);

            \draw[warc] (-12, 6) -- (-8, 6) node[midway, above, black, yshift=0pt, xshift=-1pt]{\scriptsize $n$};
            
            \draw[warc] (-12, -6) -- (-8, -6) node[midway, below, black, yshift=0pt, xshift=-2pt]{\scriptsize $m$};
            
            \draw[warc] (18, 6) -- (14, 6) node[midway, above, black, yshift=0pt, xshift=1pt]{\scriptsize $n$};
            
            \draw[warc] (18, -6) -- (14, -6) node[midway, below, black, yshift=0pt, xshift=2pt]{\scriptsize $m$};

            \draw[warc] (18, -6) to[out=east, in=east] (18, -12);
            \draw[warc] (-12, -6) to[out=west, in=west] (-12, -12);
            \draw[warc] (-12, -12) -- (18, -12);

            \draw[warc] (18, 6) to[out=east, in=east] (18, 12);
            \draw[warc] (-12, 6) to[out=west, in=west] (-12, 12);
            \draw[warc] (-12, 12) -- (18, 12);
             
            \node[vJW, scale=0.06] at (-8,6) {};
            \node[vJW, scale=0.06] at (-8,-6) {};
            \node[vJW, scale=0.06] at (14,6) {};
            \node[vJW, scale=0.06] at (14,-6) {};
        \end{scope}
    }
} \mkern-8mu \Biggr \rangle$$
The final diagram resolves to $(-1)^{n+K_l} \frac{[n+1][m+1]}{[m-K_l +1]}$ by the closure property of Lemma 7.3. \hfill$\square$

This gives us the 2-colored Jones polynomial for this case of odd $(m,n)-$colored 2-bridge links. For $m=n$, this reduces to the $n$-colored Jones polynomial of odd 2-bridge links treated by Yuasa in \cite{yuasa}. Swapping the colors $m$ and $n$ on the two components gives a second formula; we expect it to follow from the same computation, but we do not carry this out here.\

\section*{Further Discussions}

Besides just the odd 2-bridge links, there are more classes of knots which we think would be interesting to reduce using the generalized full-twist formula, as well as formulas for braids of colored strands as opposed to only twists. One family of links that we think would be worth exploring is the class of pretzel links for which we could find formulae for arbitrarily many differently colored components:
$$\mbox{
    \tikz[baseline=-.6ex, scale=.06]{
        \setsquarecoords{12}
         
        \begin{scope}

            \draw[warc] (-4, 8) to[out=north, in=north] (-20, 8);
            \draw[warc] (4, 8) to[out=north, in=north] (20, 8);
            \draw[warc] (-28, 8) to[out=north, in=north] (-44, 8);
            \draw[warc] (28, 8) to[out=north, in=north] (44, 8);

            \draw[warc] (-4, -8) to[out=south, in=south] (-20, -8);
            \draw[warc] (4, -8) to[out=south, in=south] (20, -8);
            \draw[warc] (-28, -8) to[out=south, in=south] (-44, -8);
            \draw[warc] (28, -8) to[out=south, in=south] (44, -8);

            \draw[warc] (52, 8) .. controls (52, 48.7) and (-52, 48.7) .. (-52, 8);
            \draw[warc] (52, -8) .. controls (52, -48.7) and (-52, -48.7) .. (-52, -8);

            \draw[fill=white] (56,-8) rectangle (40,8) node[midway] {$2a_1$};
            \draw[fill=white] (32,-8) rectangle (16,8) node[midway] {$2a_2$};
            \draw[fill=white] (-8,-8) rectangle (8,8) node[midway] {$\cdots$};
            \draw[fill=white] (-32,-8) rectangle (-16,8) node[midway] {$2a_{l-1}$};
            \draw[fill=white] (-56,-8) rectangle (-40,8) node[midway] {$2a_l$};

        \end{scope}
    }
}$$
\clearpage
\section*{Acknowledgements}
I would also like to give acknowledgments, firstly to Ka Ho Wong for being my advisor over the course of this paper, who has served as my introduction into knot theory and the volume conjecture, and been a passionate and willing mentor. Secondly to Wataru Yuasa, who was willing to share his workflow for creating LaTeX diagrams and whose paper has served as an inspiration for deriving generalized formulas for multiple colors.\\

AI Acknowledgment: Claude was used to revise text and create citations. The content of the paper and its methods are human-crafted.


\begin{thebibliography}{9}
\bibitem{ChenZhu} Q.~Chen and S.~Zhu, \emph{On the asymptotic expansions of various quantum invariants I: the colored Jones polynomial of twist knots at the root of unity $e^{\frac{2\pi\sqrt{-1}}{N+\frac{1}{2}}}$}, preprint, \texttt{arXiv:2307.12963 [math.GT]}.
\bibitem{habiro} K.~Habiro, \emph{On the quantum $\spl_2$ invariants of knots and integral homology spheres}, Geom. Topol. Monogr. \textbf{4} (2002), 55--68, \texttt{arXiv:math/0211044}.
\bibitem{KT} R.~M.~Kashaev and O.~Tirkkonen, \emph{Proof of the volume conjecture for torus knots}, Zap. Nauchn. Sem. POMI \textbf{269} (2000), 262--268; English transl.\ J. Math. Sci. (N.Y.) \textbf{115} (2003), no.~1, 2033--2036, \texttt{arXiv:math/9912210}.
\bibitem{kauffman} L.~H.~Kauffman and S.~Lins, \emph{Temperley--Lieb Recoupling Theory and Invariants of 3-Manifolds}, Annals of Mathematics Studies 134, Princeton University Press, 1994.
\bibitem{masbaum} G.~Masbaum, \emph{Skein-theoretical derivation of some formulas of Habiro}, Algebr. Geom. Topol. \textbf{3} (2003), 537--556, \texttt{arXiv:math/0306345}.
\bibitem{MM} H.~Murakami and J.~Murakami, \emph{The colored Jones polynomials and the simplicial volume of a knot}, Acta Math. \textbf{186} (2001), no.~1, 85--104.
\bibitem{JunMurakami} J.~Murakami, \emph{Complexified tetrahedrons, fundamental groups, and volume conjecture for double twist knots}, preprint, \texttt{arXiv:2501.00225 [math.GT]}.
\bibitem{yuasa} W.~Yuasa, \emph{The $\spl_3$ colored Jones polynomials for 2-bridge links}, preprint, \texttt{arXiv:1609.07289}.
\end{thebibliography}
\end{document}